\documentclass[a4j,12pt]{article} 
\usepackage{amsmath,amssymb}

\newcommand{\Section}[1]{%
\renewcommand{\thesection}{\S\arabic{section}}
\section{#1}
\renewcommand{\thesection}{\arabic{section}}
\setcounter{equation}{0}}

\newtheorem{lem}{Lemma}[section]
\newtheorem{rem}[lem]{Remark}
\newtheorem{thm}[lem]{Theorem}  %' 'é'©'à
\newtheorem{prop}[lem]{Proposition}
\newtheorem{cor}[lem]{Corollary}

\newcommand{\B}{{\Bbb B}}
\newcommand{\G}{{\Bbb G}}
\newcommand{\Q}{{\Bbb Q}}
\newcommand{\Z}{{\Bbb Z}}

\newcommand{\C}{{\Bbb C}}

\newcommand{\bH}{{\Bbb H}} 
\newcommand{\bP}{{\Bbb P}} 

\newcommand{\T}{{\Bbb T}}
\newcommand{\V}{{\Bbb V}}
\newcommand{\X}{{\Bbb X}}
\newcommand{\Y}{{\Bbb Y}}

\newcommand{\calC}{{\cal C}}
\newcommand{\calD}{{\cal D}}
\newcommand{\calF}{{\cal F}}
\newcommand{\calH}{{\cal H}}
\newcommand{\calO}{{\cal O}}
\newcommand{\calP}{{\cal P}}
\newcommand{\calR}{{\cal R}}
\newcommand{\calS}{{\cal S}}
\newcommand{\calU}{{\cal U}}

\newcommand{\frX}{{\frak X}}
\newcommand{\frp}{{\frak p}}

\newcommand{\lam}{{\lambda}}
\newcommand{\Lam}{{\Lambda}}
\newcommand{\ve}{{\varepsilon}}
\newcommand{\alp}{{\alpha}}  %2003
\newcommand{\vphi}{{\varphi}} %2004
\newcommand{\eps}{{\epsilon}} %2008

\newcommand{\slit}{\vspace{5mm}}
\newcommand{\mslit}{\vspace{3mm}}

\newcommand{\qed}{\hfill \hbox{\rule[-2pt]{4pt}{7pt}}}
\newcommand{\proof}{{\hspace*{0.4cm} {\it Proof}.\ \enskip}}

\newcommand{\hec}{{\calH(G,K)}}
\newcommand{\CKX}{{\calC^\infty(K \backslash X)}}
\newcommand{\SKX}{{\calS(K \backslash X)}}
\newcommand{\SXtilde}{{ \calS(\wt{X})}}
\newcommand{\SX}{{ \calS(X)}}

\newcommand{\real}{{\rm Re}}
\newcommand{\Hom}{{\rm Hom}}

\newcommand{\set}[2]{\left\{\left.#1\vphantom{#2}\:\right\vert\:#2\right\}}
\newcommand{\wt}{\widetilde}
\newcommand{\what}{\widehat}

\newcommand{\ddprod}[2]{\displaystyle{\prod_{#1}^{#2}}}
\newcommand{\dsum}[1]{\displaystyle{\sum_{#1}}}
\newcommand{\dprod}[1]{\displaystyle{\prod_{#1}}}
\newcommand{\dsqcup}[1]{\displaystyle{\bigsqcup_{#1}}}
\newcommand{\dint}[1]{\displaystyle{\int_{#1}}}

\newcommand{\rank}{{\rm rank}}
\newcommand{\abs}[1]{\left\vert{#1}\right\vert}  %in mymacro.tex

\newcommand{\Supp}{{\rm Supp}}
\newcommand{\CG}{{\calC^\infty(G)}} 
\newcommand{\DG}{{\calD(G)}} 
\newcommand{\SG}{{\calS(G)}} 

\newcommand{\pair}[2]{{\langle}{#1},{#2}{\rangle}}

\newcommand{\twomatrix}[4]{\left(\begin{array}{cc}
                           {#1} & {#2}\\
                           {#3} & {#4}
                          \end{array}\right)}
\newcommand{\twovector}[2]{\left(\begin{array}{c}{#1}\\{#2}\end{array}\right)}

\newcommand{\twomatrixminus}[4]{\begin{array}{cc}
                           {#1}  & {#2}\\
                           {#3} &  {#4}
                          \end{array}}

\newcommand{\gyaddots}{%
  \setlength{\unitlength}{1mm}
   \begin{picture}(5,3.5)(-2,-1.5)
   \put(0,0){$\cdot$}
   \put(2,1.5){$\cdot$}
   \put(-2,-1.5){$\cdot$}
   \end{picture}}
\newcommand{\gyakuddots}{\smash{\lower0.3ex\hbox{\gyaddots}}}
\begin{document}
\title{{\large Spherical functions on $p$-adic homogeneous spaces}}
\author{Yumiko Hironaka}
\date{}
\maketitle

%\hfill{\today}

\renewcommand{\thefootnote}{{}}

\footnotetext{
\hspace*{-2.8mm}2000 {\it Mathematics Subject Classification.} Primary 11F85; Secondly 11E95, 11F70, 22E50. 


This research was partially supported by Grant-in-Aid for Scientific Research (C):20540029.

}

%\begin{center}
%\begin{minipage}{14cm}
%{\small {\bf Abstract.} 

%Let ${\Bbb G}$ be a connected reductive linear algebraic group and ${\Bbb X}$ a ${\Bbb G}$-homogeneous affine algebraic variety both defined over a ${\frak p}$-adic field $k$, where we assume a minimal $k$-parabolic subgroup of ${\Bbb G}$ acts with open orbit. We are interested in spherical functions on $X = \X(k)$.
%, which are nonzero common eigenfunctions on $X$ with respect to the action of Hecke algebra ${\cal H}(G, K)$ and $K$-invariant, where $K$ is a compact open subgroup of $G = \G(k)$.
%
%In the present papaer, we give a unified method to obtain functional equations of spherical functions on $X$ under the condition (AF) in Introduction, and explain functional equations are reduced to those of $p$-adic local zeta functions of small prehomogeneous vector spaces of limited type. 
%}
%\end{minipage}
%\end{center}

\setcounter{section}{0}

\vspace{5mm}
%\input{intro}

%%%%% intro %%%%

\noindent
{\Large {\bf Introduction}}

\medskip
\noindent
Let $\G$ be a reductive linear algebraic group defined over $k$, 
and $\X$ be an affine algebraic variety defined over $k$ which is $\G$-homogeneous, 
where and henceforth $k$ stands for a non-archimedian local field of characteristic $0$.
The Hecke algebra $\hec$ of $G$ with respect to $K$ acts by convolution product on the space of $\CKX$ of $K$-invariant $\C$-valued functions on $X$, where $K$ is a maximal compact open subgroup of $G = \G(k)$ and $X = \X(k)$.

A nonzero function in $\CKX$ is called {\it a spherical function on $X$} if it is a common $\hec$-eigen function.

\slit
Spherical functions on homogeneous spaces comprise an interesting topic to investigate and a basic tool to study harmonic analysis on $G$-space $X$. 
They have been studied also as spherical vectors of distinguished models, Shalika functions and Whittaker-Shintani functions, and are closely related to theory of automorphic forms and representation theory. 
When $\G$ and $\X$ are defined over $\Q$, spherical functions appear in local factors of global objects, e.g. Rankin-Selberg convolutions and Eisenstein series (e.g. [CS], [Fl], [HS3], [Jac], [KMS], [Sf2]). 

The theory of spherical functions also has applications to classical number theory. For example when $X$ is the space of symmetric forms, alternating forms or hermitian forms, spherical functions can be considered as generating functions of local densities, and have been applied to obtain their explicit formulas (cf. [HS1], [HS2], [H1]-[H4]).

To obtain explicit expressions of spherical functions is one of basic problems. For the group cases, it has been done by I.~G.~Macdonald and afterwards by W.~Casselman by a representation theoretical method (cf. [Ma], [Cas]). There are some results on homogeneous space cases mainly for the case that the space of spherical functions attached to each Satake parameter is of dimension one (e.g. [CS], [KMS], [Of]). 

\newpage
In this paper, following the preliminaries in \S 1, we give a general expression of spherical functions on $X$ of dimension not necessary one based on the data of the group $G$ and functional equations of spherical functions in \S 2.
Then we show a unified method to obtain functional equations of spherical functions on $X$, and explain that functional equations are reduced to those of $p$-adic local zeta functions of {\it small} prehomogeneous vector spaces in \S 3. 
These are improvements of some results in \cite{JMSJ} and \cite{Hamb}. We devote \S 4 to examples.
For general references for algebraic groups, one may refer to \cite{Bo} and \cite{PR}.
\slit

\bigskip
\Section{}%Preliminaries}

\mslit
{\bf 1.1.} 
Let $\bH$ be a connected linear algebraic group and $\Y$ an affine algebraic variety on which $\bH$ acts, where everything is assumed to be defined over $k$.
We denote by $\frX(\bH)$ the group of $k$-rational characters of $\bH$, which is a free abelian group of finite rank. We set $\frX_0(\bH)$ for the subgroup consisting of characters corresponding to some relative $\bH$-invariants on $\Y$, where a rational function $f$ on $\Y$ defined over $k$ is called {\it relative $\bH$-invariant} if it satisfies, for some $\psi \in \frX(\bH)$,  
$$
f(g \cdot y) = \psi(g) f(y), \quad g \in \bH.
$$
We say a set $\set {f_i(y)} {1 \leq i \leq n}$ is {\it basic}, if the corresponding characters form a basis for $\frX_0(\bH)$; 
then every relative $\bH$-invariant on $\Y$ has a following form
$$
c \cdot \ddprod{i=1}{n}\, f_i(y) ^{e_i}, \qquad c \in k^\times, \; e_i \in \Z.
$$
We consider the following conditions for $(\bH, \Y)$.

\mslit
(A1') $\Y$ has a Zariski open $\bH$-orbit. 

\mslit
(A1) $\Y$ has only a finite number of $\bH$-orbits.

\mslit
(A2) A basic set of relative $\bH$-invariants on $\Y$ can be taken by regular functions on $\Y$. 

\mslit
(A3) For $y \in \Y$ not contained in open orbits, there exists some $\psi$ in  $\frX_0(\bH)$ whose restriction to the identity component of the stabilizer $\bH_y$ is not trivial.  

\mslit
(A4) The rank of $\frX_0(\bH)$ coincides with that of $\frX(\bH)$.

\begin{rem}
{\rm Assume that $\Y$ is a homogeneous space of a connected reductive linear algebraic group $\G$ (like as \S 1. 2). Then $\Y$ is irreducible and there is at most one Zarisky open orbit, and (A1) implies (A1').
The condition (A1) is satisfied if $\Y$ is a spherical homogeneous space of $\G$ and $\bH$ is a minimal parabolic subgroup of $\G$, and symmetric spaces are spherical, especially the spaces of type $\G/\G^\theta$ ($\G^\theta$ is the set of fixed points of an involution $\theta$ on $\G$) are spherical(cf. \cite{Sf4}).
As for (A2), we note here that for the case of prehomogeneous vector spaces, basic relative invariants can be chosen as polynomial functions (\cite{Sf1}-Lemma 1.3). The condition (A3) assures us a good condition for distributions on $\Y(k)$, which we need when we consider functional equations of spherical functions.   
}%
\end{rem}

\bigskip
\noindent
{\bf 1.2.} 
Hereafter, let $\G$ be a connected reductive linear algebraic group and $\X$ be an affine algebraic variety which is $\G$-homogeneous, where everything is assumed to be defined over $k$. For an algebraic set, we use the same ordinary letter for the set of $k$-rational points, e.g. $G = \G(k), \; X = \X(k)$.
Let $K$ be a maximal compact open subgroup of $G$, and $\B$ a minimal parabolic subgroup of $\G$ defined over $k$ satisfying $G = KB = BK$. The group $\B$ is not necessarily a Borel subgroup. 
We denote by $\abs{\; }$ the absolute value on $k$ normalized by $\abs{\pi} = q^{-1}$, where $\pi$ is a prime element of $k$ and $q$ is the cardinal number of the residue class field of $k$, we understand $\abs{0} = 0$ for simplicity.

\mslit
Assume that $(\B, \X)$ satisfies (A1') and (A2), and let $\set{f_i(x)}{1 \leq i \leq n}$ be a basic set of regular relative $\B$-invariants, and $\psi_i \in \frX_0(\B)$ the corresnponding character to $f_i(x)$, and $n = \rank(\frX_0(\B))$. The open $\B$-orbit $\X^{op}$ is decomposed into a finite number of open $B$-orbits over $k$ (cf. \cite{Serre}-III-4.4), which we write 
$$
\X^{op}(k) = \dsqcup{u \in J(X)}\, X_u.
$$  

For $x \in X$, $s \in \C^n$ and $u \in J(X)$, we define
\begin{eqnarray} \label{def:sph}
\omega(x; s) = \dint{K}\, \abs{f(k\cdot x)}^s dk,\quad
\omega_u(x; s) = \dint{K}\, \abs{f(k\cdot x)}_u^s dk,
\end{eqnarray}
where 
$dk$ is the Haar measure on $K$ normalized  by $\dint{K} dk = 1$, and 
\begin{eqnarray} \label{def:f^s}
\abs{f(x)}^{s} = \ddprod{i=1}{n}\, \abs{f_i(x)}^{s_i},
\quad
\abs{f(x)}_u^{s} =\left\{ \begin{array}{ll}
   \abs{f(x)}^{s}  & \mbox{if } x \in X_u,\\
   0 & \mbox{otherwise}.
\end{array}
 \right. \nonumber
\end{eqnarray}
We set 
$$
\abs{\psi(p)}^s = \prod_{i=1}^n \abs{\psi_i(p)}^{s_i}.
$$
By the following proposition, we see $\omega(x;s)$ and $\omega_u(x; s)$ are spherical functions on $X$, where we give also the 'eigenvalues' for them.  

\begin{prop} \label{prop in S1}
The integrals in (\ref{def:sph}) are absolutely convergent if $\real(s_i) \geq 0, \; 1 \leq i \leq n$, analytically continued to rational functions of $q^{s_1}, \ldots, q^{s_n}$, and become $\hec$-common eigen functions. In particular $\omega_u(x; s), \; u \in J(X)$, are spherical functions on $X$ and linearly independent for generic $s$. More precisely, 
for each $\phi \in \hec$, one has $(\phi * \omega(\;;s))(x) = \lam_s(\phi) \omega(x;s)$ and $(\phi * \omega_u(\;;s))(x) = \lam_s(\phi) \omega_u(x;s)$ with
\begin{eqnarray*}
\lam_s(\phi) = \dint{B} \phi(p) \abs{\psi(p)}^{-s}\delta(p)dp = \dint{G} \phi(g) \abs{\psi(p(g))}^{-s} \delta(p(g)) dg, 
\end{eqnarray*}
where $dp$ is the left invariant Haar measure on $B$ normalized by $\dint{K \cap B} dp = 1$ and $p(g) \in B$ for which $p(g)^{-1}g \in K$.
\end{prop}

\medskip
{For a proof} we refer to \cite{JMSJ}-Proposition 1.1. 
To make it sure we note here the action of $\hec$ on $\omega_u(x;s)$ : for $\phi \in \hec$ and $x \in X$, 
\begin{eqnarray*}
(\phi * \omega_u(\; ; s))(x) &=& 
\dint{G}\phi(g)\dint{K}\abs{f(kg^{-1}\cdot x)}_u^{s}dkdg \\ 
&=&
\dint{K} \dint{G} \phi(gk) \abs{f(g^{-1}\cdot x)}_u^{s} dg dk
=
\dint{G} \phi(g) \abs{ f(g^{-1}\cdot x)}_u^{s} dg\\
&=&
 \dint{K} \dint{B} \phi(kp) \abs{f(p^{-1}k^{-1}\cdot x)}_u^{s} d_rp dk\\
&=&
\dint{K} \dint{B} \phi(p) \abs{\psi(p)}^{-s} \abs{f(k^{-1}\cdot x)}_u^{s} d_rp dk\\
& =&
\dint{B} \phi(p) \abs{\psi(p)}^{-s} \delta(p)dp \cdot \omega_u(x; s).
\end{eqnarray*}
\qed

\bigskip
\begin{rem} \label{modified sph}
{\rm When we assume also (A4), we can determine $\ve_0 \in \Q^n$ by 
\begin{eqnarray*}
\abs{\psi(p)}^{\ve_0} = \delta^{\frac12}(p), \; p \in B,
\end{eqnarray*}
and it is better to modify the definition of spherical functions as follows:
\begin{eqnarray} \label{modified sph}
\wt{\omega}_u(x; s) = \dint{K}\, \abs{f(k\cdot x)}_u^{s+\ve_0}dk.
\end{eqnarray}
Then, instead of Proposition~\ref{prop in S1}, we have 
\begin{eqnarray*}
\left( \phi * \wt{\omega}(\; ; s) \right)(x) = \wt{\lam}_s(\phi) \wt{\omega}(x; s), 
\quad \left( \phi * \wt{\omega}_u(\; ;  s) \right)(x) = \wt{\lam}_s(\phi) \wt{\omega}_u(x; s), 
\end{eqnarray*}
where
\begin{eqnarray*}
\wt{\lam}_s(\phi) = \dint{B} \phi(p) \abs{\psi(p)}^{-s+\ve_0}dp
= \dint{G}\phi(g) \abs{\psi(p(g))}^{-s+\ve_0} dg, \quad  \phi \in \hec.
\end{eqnarray*}
} % end \rm
\end{rem}  

\medskip
\begin{rem} \label{rem sph with char}
{\rm 
The value $f_i(x)\mod\psi_i(B)$ is constant in $k^\times/\psi_i(B)$ on each open $B$-orbit, which we call {\it the signature of } $f_i$. If we can parametrize open $B$-orbits in $X^{op}$ by the signatures of $\set{f_i(x)}{1 \leq i \leq n}$, then $J(X)$ can be naturally identified with a subset of the finite abelian group 
\begin{eqnarray} \label{sign gr}
\left(  k^\times \right)^n \big{/} \ddprod{i=1}{n}\, \psi_i(B).
\end{eqnarray}
Such cases often occur, and then, it is natural to consider spherical functions with character as follows. 
Assume $\calU$ is $J(X)$ or its suitable subset which is canonically identified with a subgroup of (\ref{sign gr}). Taking a character $\chi$ of $\calU$, we set
\begin{eqnarray} 
&&
\omega(x; \chi; s) = 
\dint{K}\, \chi(f(k\cdot x)) \abs{f(k\cdot x)}^s dk =
\dsum{u \in \calU}\, \chi(u) \omega_u(x; s), \label{char-sph} \\
&&
\wt{\omega}(x; \chi; s) = \dint{K}\, \chi(f(k\cdot x)) \abs{f(k\cdot x)}^{s+\ve_0} dk =
\dsum{u \in \calU}\, \chi(u) \wt{\omega}_u(x; s), \label{char-mod-sph}
\end{eqnarray}
the latter can be considered only when (A4) is satisfied. 
} %end \rm
\end{rem}

\bigskip
Let $W$ be the relative Weyl group of $\G$ with respect to $\T$, where $\T$ is a maximal $k$-split torus contained in $\B$. 
The group $W$ acts on $\frX(\B)$ as $(\sigma\xi)(b) = \xi(n_\sigma^{-1}bn_\sigma)$ by taking a representative $n_\sigma \in Z_G(T)$ of $\sigma \in W$, hence it acts on $s \in \C^n$ through the identification $\C^n \cong \frX_0(\B) \otimes_\Z \C \subset \frX(\B) \otimes_\Z \C$.

\vspace{2cm}
\Section{}  %\S 2
In this section, we will give a general expression for spherical functions based on the data of the group $G$ and functional equations of spherical functions.

We follow the notation in \S 1.2, and take $K$ as a special, good, maximal compact subgroup in the sense of Bruhat and Tits (cf. \cite{Cas}-\S 3.5), and Iwahori subgroup $U$ of $K$ compatible with $B$. 

\bigskip
\noindent
{\bf 2.1.} 
In this subsection, we prepare some results from representation theory (cf. \cite{JMSJ}-\S 1, \cite{Cas}).

We denote by $\SG$ the Schwartz-Bruhat space on $G$, namely the space of locally constant compactly supported functions on $G$, and set $\DG = \Hom_\C(\SG, \C)$, the space of distributions on $G$, and the  pairing on $\DG \times \SG$
\begin{eqnarray*}
\pair{T}{\phi} = \pair{T}{\phi}_{\calD \times \calS} = T(\phi), \quad (T\in \DG, \; \phi \in \SG).
\end{eqnarray*}
Then the space $\CG$ of locally constant functions can be regarded as a subspace of $\DG$ by
\begin{eqnarray*}
\pair{\psi}{\phi} = \dint{G}\, \psi(g)\phi(g)dg, \quad (\psi \in \CG, \; \phi \in \SG).
\end{eqnarray*}
We regard $\CG$ as a two-sided $G$-module and $\SG$ as a submodule by  
\begin{eqnarray*}
g \cdot \psi(x) = \psi(xg), \qquad
\psi^g(x) = \psi(gx), \quad (\psi \in \CG, \; g,x \in G).
\end{eqnarray*}
Then,  $\DG$ becomes also a two-sided $G$-module by the dual action:
\begin{eqnarray*}
\pair{g \cdot T}{\phi} = \pair{T}{g^{-1} \cdot \phi}, \quad 
\pair{T^g}{\phi} =  \pair{T}{\phi^{g^{-1}}}, \quad
(T \in \DG, \; \phi \in \SG, \; g \in G).
\end{eqnarray*}

For a subspace $\Gamma$ of $\DG$ and a subgroup $H$ of $G$, we denote by $\Gamma^H$ the set of left $H$-invariant elements in $\Gamma$. 

\bigskip
Let $\chi$ be an unramified regular character of the centralizer $Z_G(T)$, i.e. $\chi\vert_{Z_G(T) \cap K} \equiv 1$, and $\sigma \chi = \chi$ implies $\sigma = 1$ for $\sigma \in W$, which is canonically extended to be a character of $B$.
We recall the induced representation (principal series representation) of $G$:
\begin{eqnarray}
I(\chi) 
&=& {\rm Ind}_{B}^{G} (\chi) = \set{\phi \in \CG} 
{\phi(pg) = \chi\delta^{\frac12}(p) \phi(g) \quad(p \in B, \; g \in G)}\nonumber\\
&=&
\set{\phi \in \CG}{\phi^p = \chi\delta^{\frac12}(p)\phi, \quad (p \in B)},
\end{eqnarray}
which is a left $G$-submodule of $\CG$.
Then we have a left $G$-equivariant surjection 
\begin{eqnarray*} \label{P_chi}
&&\calP_\chi: \SG \longrightarrow I(\chi), \\
&&
\calP_\chi(\phi)(x) = \dint{B}\, \chi^{-1}\delta^\frac12 (p) \phi(px)dp, \quad (\phi \in \SG, \; x \in G). \nonumber
\end{eqnarray*}
We set $\vphi_{K, \chi} = \calP_\chi(ch_K)$, where $ch_K$ is the characteristic function of $K$.

%\medskip
The map $\calP_{\chi^{-1}}$ induces a left $G$-equivariant injection ${\calP_{\chi^{-1}}}^*$ from $I(\chi^{-1})^* = \Hom_\C(I(\chi^{-1}), \C)$ to $\DG$ determined by
\begin{eqnarray*}
&&
\pair{{\calP_{\chi^{-1}}}^* (T)}{\phi} = T(\calP_{\chi^{-1}}(\phi)),   \quad (T \in I(\chi^{-1})^*, \; \phi \in \SG), 
\end{eqnarray*}
and we obtain the following (\cite{JMSJ}-Lemma 1.2, Corollary 1.3) \label{JMSJ-Lem.1.2}).

\begin{prop}  \label{I and DG}
By the dual map ${\calP_{\chi^{-1}}}^*$ of $\calP_{\chi^{-1}}$, one has a left $G$-isomorphism
\begin{eqnarray*}
I(\chi^{-1})^* \cong \calD(G)_\chi := \set{T \in \DG}{T^p = \chi\delta^{\frac12}(p)T \quad (p \in B)}.
\end{eqnarray*}
Further, by this isomorphism, $I(\chi)$ and $I(\chi^{-1})$ can be understood as the smooth dual of the each other with pairing 
%\begin{eqnarray*} \label{smooth dual}
%\CG \cap I(\chi^{-1})^* \cong \CG \cap \DG_\chi = I(\chi),
%\end{eqnarray*}
\begin{eqnarray*}
<f_1, f_2> = \dint{K}\, f_1(k)f_2(k)dk, \qquad (f_1 \in I(\chi), \; f_2 \in I(\chi^{-1})).
\end{eqnarray*}
\end{prop}

\medskip
Indeed we calculate the pairing on $I(\chi) \times I(\chi^{-1})$ in the following: for $(f_1, f_2) \in I(\chi) \times I(\chi^{-1})$ and $\phi \in \SG$ such that $\calP_{\chi^{-1}}(\phi) = f_2$,
\begin{eqnarray*}
\pair{f_1}{f_2} &=& \pair{({\calP_{\chi^{-1}}}^*)^{-1}(f_1)}{\calP_{\chi^{-1}}(\phi)}_{I(\chi^{-1})^* \times I(\chi^{-1})}\\
&=&
\pair{f_1}{\phi}_{\calD \times \calS} = \dint{G}f_1(g)\phi(g) dg = \dint{K} \dint{B}f_1(pk)\phi(pk)dp dk\\
&=&
\dint{K} f_1(k) \dint{B}\chi\delta^{\frac12}(p)\phi(pk)dp dk = \dint{K}f_1(k)\calP_{\chi^{-1}}(\phi)(k)dk\\
&=&
\dint{K}f_1(k)f_2(k)dk.
\end{eqnarray*}

\bigskip
For $\sigma \in W$, there is a unique left $G$-equivariant map satisfying 
\begin{eqnarray}
T_\sigma^\chi : I(\chi) \longrightarrow I(\sigma\chi),  \qquad
T_\sigma^\chi(\vphi_{K,\chi}) = c_\sigma(\chi) \vphi_{K,\sigma\chi},
\end{eqnarray}
where
\begin{eqnarray}
c_\sigma(\chi) &=& \dprod{\alp \in \Sigma^+, \; \sigma(\alp) < 0}\, c_\alp(\chi), \nonumber \\
c_\alp(\chi) &=& \dfrac{(1 - q_{\frac{\alp}{2}}^{-\frac12}q_\alp^{-1}\chi(a_\alp))
(1 + q_{\frac{\alp}{2}}^{-\frac12}\chi(a_\alp))}
{1 - \chi(a_\alp)^2} \nonumber \\
&(=&
\dfrac{1 - q^{-1}\chi(a_\alp)}{1 - \chi(a_\alp)} \quad \mbox{if $G$ is split}). \nonumber
\end{eqnarray}
Here $\Sigma^+$ is the set of positive roots $\G$ with respect to $\T$ and $\B$, and for the definition of $a_\alp \in T$ and numbers $q_\alp, \; q_{\frac12 \alp} \; (\alp \in \Sigma)$, see \cite{Cas}.
It is known that
\begin{eqnarray} \label{isom}
T_\sigma^\chi \mbox{ is an isomorphism if and only if } c_\sigma(\chi)c_{\sigma^{-1}}(\sigma(\chi)) \ne 0.
\end{eqnarray}

For a compact open subgroup $V$ of $G$, we define an operator $\calP_V$ on $\DG$ by
\begin{eqnarray*}
\pair{\calP_V(T)}{\phi} = \dint{V}\pair{u\cdot T}{\phi} du = \dint{V}\pair{T}{u^{-1}\cdot \phi}du, \quad (T \in \DG, \; \phi \in \SG),
\end{eqnarray*}
where $du$ is the Haar measure on $V$ normalized by $\dint{V} du = 1$.

As the adjoint $G$-morphism of $T_{\sigma^{-1}}^{\sigma \chi^{-1}}$, we have (under the identification through ${\calP_{\chi}}^*$ and ${\calP_{\sigma\chi^{-1}}}^*$ by Proposition~\ref{I and DG})
$$
\left(T_{\sigma^{-1}}^{\sigma \chi^{-1}} \right)^* : \DG_{\chi} = I(\chi^{-1})^* \longrightarrow \DG_{\sigma\chi} = I(\sigma \chi^{-1})^*  .
$$

Then we see the following(cf. \cite{JMSJ}-Propposition1.6, Proposition 1.7).

\begin{prop} \label{extension of T}
Assume $c_\sigma(\chi) c_\sigma(\chi^{-1})c_{\sigma^{-1}}(\sigma \chi)c_{\sigma^{-1}}(\sigma \chi^{-1}) \ne 0$. Then
$$
\wt{T}_\sigma^\chi = \frac{c_\sigma(\chi)}{c_{\sigma^{-1}}(\sigma \chi^{-1})} 
\left( T_{\sigma^{-1}}^{\sigma \chi^{-1}}  \right)^* : \DG_\chi \longrightarrow \DG_{\sigma \chi}
$$
is an extension of the $G$-isomorphism $T_\sigma^{\chi} : I(\chi) \longrightarrow I(\sigma \chi)$.
Further, for a compact open subgroup $V$ of $G$, one has
$$
\calP_V \circ \wt{T}_\sigma^{\chi} = T_\sigma^\chi \circ \calP_V.
$$
\end{prop}

\bigskip
We recall Casselman basis $\set{f_{\sigma,\chi}}{\sigma \in W}$ for $I(\chi)^U$, which satisfies the following(cf. \cite{Cas})
\begin{eqnarray*}
&&
T_\sigma^\chi(f_{\tau,s})(1) = \delta_{\sigma, \tau}, \\
&&
\calP_K(f_{\sigma,\chi})(1) = 
\frac{\gamma(\sigma(\chi))}{Q \cdot c_\sigma(\chi)},
\end{eqnarray*}
where $\delta_{\sigma, \tau}$ is the Kronecker delta, and
\begin{eqnarray*}
\gamma(\chi) = \dprod{\alp \in \Sigma^+}\, c_\alp(\chi),\quad
Q = \dsum{\sigma \in W}\, [U\sigma U : U]^{-1}.
\end{eqnarray*}

\bigskip
Let us recall our situation in \S 1. 
Since $\G$ is reductive and $\X$ is affine, $\bH = \G_{x_0}, \; (x_0 \in X^{op})$ is reductive also(cf. \cite{Sp}-Satz 3.3).
Then there is a $G$-invariant measure on $G/H$, since $G$ and $H$ are unimodular.
Since $BH/H$ is open in $G/H$ and isomorphic to $B/B_0$ with $B_0 = B \cap H$, there is an invariant measure on $B/B_0$, so the modulus character of $B_0$ coincides with $\delta\vert_{B_0}$.

Now we set
$$
I(\chi, BH) = \set{\phi \in I(\chi)}{Supp(\phi) \subset BH}.
$$
The next lemma is based on an idea of O.~Offen used in \cite{Of}. It will play a key role to restrict the summation with respect $W$ to a certain subgroup $W_0$ in \S 2.2. 

\begin{lem} \label{lem offen-type}
If there is a nonzero left $H$-invariant distribution in $I(\chi)^*$ which is not identically zero on $I(\chi, BH)$,
then $\chi = \delta^\frac12$ on $B_0$. 
\end{lem}

\proof
Assume a distribution $\Lambda \in I(\chi)^*$ satisfies the condition as above.
The space $I(\chi, BH)$ can be identified (by the restriction) with 
\begin{eqnarray*}
Ind_{B_0}^{H}(\chi\vert_{B_0})&=&
\set{f \in \calC^\infty(H)}{
f(p_0h) = \chi\delta^{\frac12}(p_0) f(h) \quad(p_0 \in B_0, \; h \in H)},
\end{eqnarray*}
on which there is a left $H$-invariant surjection $\calP_{H,\chi}$ from $\calS(H)$ given by
\begin{eqnarray*}
\calP_{H,\chi}(\vphi)(h) = \dint{B_0}\chi^{-1}\delta^{\frac12}(p_0) \vphi(p_0h)dp_0,
\end{eqnarray*}
where $dp_0$ is a left invariant Haar measure on $B_0$.
Then we have a nonzero $H$-invariant distribution $T$ on $H$ determined by
\begin{eqnarray*}
\pair{T}{\vphi}_{\calD \times \calS} = \pair{\Lambda}{\calP_{H,\chi}(\vphi)},
\end{eqnarray*}
thus we have a left invariant measure on $H$, which becomes also right invariant since $H$ is unimodular.
On the other hand, since we have for $p \in B_0$
\begin{eqnarray*}
\pair{T^{p^{-1}}}{\vphi}_{\calD \times \calS} &=& \pair{T}{\vphi^{p}}_{\calD \times \calS}
= \pair{\Lam}{\calP_{H, \chi}(\vphi^p)} = \pair{\Lam}{\chi\delta^{-\frac12}(p)\calP_{H, \chi}(\vphi)}\\
&=&
\chi\delta^{-\frac12}(p) \pair{T}{\vphi}_{\calD \times \calS},
\end{eqnarray*}
we obtain $\chi = \delta^{\frac12}$ on $B_0$.
\qed

\vspace{.7cm}
\noindent
{\bf 2.2.} 
Take $x_0 \in X^{op}$ and set $\bH = \G_{x_0}$ and $\calU = \set{\nu \in J(X)}{G \cdot x_0 \cap X_\nu \ne \emptyset}$. 
We define the subgroup $W_0$ of W by
$$
W_0 = \set{\sigma \in W}{\sigma(\abs{\psi}^s) \equiv 1 \mbox{ and } \sigma(\delta) = \delta \mbox{ on }{B\cap H}},
$$
where $s \in \C^n$ is considered as a variable.
Though we do not assume the condition (A4), $\sigma(\abs{\psi}^s)$ is contained in $\frX_0(\B)\otimes_\Z \C\,$ if $\sigma \in W_0$, and in this case $\abs{f(x)}_u^{\sigma(s)}$ and $\abs{\psi}^{\sigma(s)} = \sigma(\abs{\psi}^s)$ are well defined. For $\sigma \in W_0$ we define $\ve_\sigma \in \Q^n$ by
\begin{eqnarray} \label{ve-sigma}
\abs{\psi}^{2\ve_\sigma} = \delta\sigma(\delta^{-1}).
\end{eqnarray}
For $s \in \C^n$, let $\chi = \chi_s$ be the charcter of $B$ given by
\begin{eqnarray} \label{chi_s}
\chi = \abs{\psi}^s\delta^{-\frac12}, \enskip i.e., \enskip \chi(p) = \abs{\psi(p)}^s \delta^{-\frac12}(p), \quad p \in B. 
\end{eqnarray}
Then we have
\begin{eqnarray*} \label{sigma-chi}
\sigma(\chi) = \abs{\psi}^{\sigma(s)+\ve_\sigma} \delta^{-\frac12}, \quad \sigma \in W_0.
\end{eqnarray*}
For each $\nu \in \calU$, we set (through the analytic continuation for general $s \in \C^n$)
\begin{eqnarray}
\Psi_\nu(x, s; g) &=& \abs{f(g\cdot x)}_\nu^s \in \calD(G)_\chi \cong I(\chi^{-1})^*,
\label{def Psi_nu}\\
\wt{\Psi}_\nu(x, s; g) &=& \calP_U(\Psi_\nu(x,s; ))(g) = \dint{U}\abs{f(gu\cdot x)}_\nu^s du \in I(\chi)^U, \nonumber
\end{eqnarray}
where we note that
\begin{eqnarray} \label{G-action on Psi}
\Psi_\nu(g_1\cdot x, s; g) = g_1 \cdot \Psi_\nu(x, s; \;)(g) = \Psi_\nu(x, s; gg_1), \quad (g, g_1 \in G).
\end{eqnarray}

\bigskip
The condition (A3) is crucial for the next lemma (cf. \cite{JMSJ}-Lemma~1.8).

\begin{lem} \label{lemma 1.8}
Assume (A1), (A2) and (A3) for $(\B, \X)$. Then for each $x \in G\cdot x_0$ and generic $s$, the set $\set{\Psi_u(x, \sigma(s)+\ve_\sigma; g)}{u \in \calU}$ forms a basis for $\DG_{\sigma\chi}^{G_x}$ for any $\sigma \in W_0$.
Here, 'generic' means to avoid a finite number of linear relations of type $\sum_{i=1}^n m_i s_i -\alp \in (\frac{2\pi\sqrt{-1}}{\log q})$ with $m_i \in \Z, \; \alp \in \C$.
\end{lem}

\bigskip
In the following, we say $s$ is {\it generic} if $s$ is generic in the sense of Lemma~\ref{lemma 1.8}, $s$ is neither a pole nor a zero of $\omega_u(x; s) \; (u \in \calU)$, $\chi = \chi_s$ is regular, and\\ $c_\sigma(\chi)c_\sigma(\chi^{-1}) c_\sigma(\sigma^{-1}\chi) c_{\sigma}(\sigma^{-1}\chi^{-1}) \ne 0$ for every $\sigma \in W_0$ (cf. Proposition~\ref{extension of T}).

We set
\begin{eqnarray}
&& \label{good elements}
\calR = \set{x \in G \cdot x_0 \cap X^{op}}{U\cdot x \subset B\cdot x_{0}},\\
&&
\calR^+ = \set{x \in \calR}{\abs{f(u\cdot x)}^s = \abs{f(x)}^s, \;(u \in U)}. \nonumber
\end{eqnarray}

Our main theorem in this section is the following, which is a refinement of \cite{JMSJ}-Proposition~1.9, where we assumed the condition (A4).  

\begin{thm}  \label{Th1}
Assume (A1), (A2) and (A3) for $(\B, \X)$ and $s$ is generic.
For $x \in \calR$, one has 
$$
\big{(} \omega_\nu(x; s) \big{)}_{\nu \in \calU} = 
\frac{1}{Q} \dsum{\sigma \in W_0}\, \gamma(\sigma(\chi)) \cdot B_\sigma(\chi) \cdot 
\left( \dint{U}\, \abs{f(u\cdot x)}_\nu^{\sigma(s)+\ve_\sigma}du \right)_{\nu \in \calU}.
$$
Moreover, if $x \in \calR^+$, one has
$$
\big{(} \omega_\nu(x; s) \big{)}_{\nu \in \calU} = 
\frac{1}{Q} \dsum{\sigma \in W_0}\, \gamma(\sigma(\chi)) \cdot B_\sigma(\chi) \cdot 
\left( \abs{f(x)}_\nu^{\sigma(s)+\ve_\sigma} \right)_{\nu \in \calU}.
$$
Here the constant $Q$ and the rational function $\gamma(\chi)$ of $q^{s_1},\ldots, q^{s_n}$ are determined by the group $G$ as in \S 2.1, and the matrix $B_\sigma(\chi)$ is determined by
the functional equation
$$
\big{(} \omega_\nu(x; s) \big{)}_{\nu \in \calU} = 
B_\sigma(\chi) \big{(} \omega_\nu(x; \sigma(s)+\ve_\sigma) \big{)}_{\nu \in \calU}.
$$ 
\end{thm}

\medskip
We give here an outline of a proof. 
By definition of $\Psi_\nu$ and $\wt{\Psi}_\nu$, we have
\begin{eqnarray}
\omega_\nu(x; s) = \dint{K}\, \Psi_\nu(x, s; k) dk = \dint{K} \wt{\Psi}_\nu(x, s; k) dk
= \calP_K (\wt{\Psi}_\nu(x, s; \; )) (1),
\end{eqnarray}
and we may write by using Casselman basis
\begin{eqnarray*}
\wt{\Psi}_\nu(x; s; \; ) &=& \dsum{\sigma \in W}\, a_{\nu, \sigma}(x;s) f_{\sigma,\chi},
\end{eqnarray*}
where
\begin{eqnarray*}
a_{\nu,\sigma}(x; s)&=& 
T_\sigma^\chi(\wt{\Psi}_\nu(x,s;\; ))(1) = (\calP_U \circ \wt{T}_\sigma^\chi)(\Psi_\nu(x, s ;\; ))(1).
\end{eqnarray*}
Now we set $\Lam = \wt{T}_\sigma^\chi(\Psi_\nu(x_0, s ;\; )) \in I(\sigma\chi^{-1})^*$, which is left $H$-invariant, and $x = g_1\cdot x_0$. Then we have 
(cf. (\ref{G-action on Psi}))
\begin{eqnarray*}
a_{\nu,\sigma}(x; s)&=& 
(\calP_U \circ \wt{T}_\sigma^\chi)(g_1\cdot \Psi_\nu(x_0, s ;\; ))(1)\\
&=&
(\calP_U \circ (g_1\cdot \Lambda))(1).
\end{eqnarray*}
Since $\calP_U(g_1 \cdot \Lambda)$ is regarded as an element of $I(\sigma\chi)^U$  by Proposition~\ref{I and DG},  
taking $\vphi_U \in I(\sigma\chi^{-1})$ as supported by $BU$ and $\vphi_U(u) = 1$ for $u \in U$, we can continue  
\begin{eqnarray*}
a_{\nu,\sigma}(x; s) &=&  \dint{K}\, \calP_U \circ (g_1\cdot \Lambda)(k)\vphi_U(k)dk\\
&=&
\pair{\calP_U \circ (g_1\cdot \Lambda)}{\vphi_U} 
= \dint{U}\pair{g_1\cdot \Lam}{u \cdot \vphi_U}du\\
&=&
\pair{g_1 \cdot \Lam}{\vphi_U} = \pair{\Lam}{g_1^{-1}\cdot \vphi_U}, 
\end{eqnarray*}
where $\Supp(g_1^{-1}\cdot \vphi_U) = BUg_1$. 

If $x \in \calR$, then $\Supp(g_1^{-1}\cdot \vphi_U) \subset BH$, and
$a_{\nu, \sigma}(x; s) = 0$ unless $\sigma\chi^{-1} = \delta^{\frac12}$ on $B_0$ by Lemma~\ref{lem offen-type}, i.e., $a_{\nu, \sigma}(x; s) = 0$ unless $\sigma \in W_0$, by our choice of $\chi$ and $W_0$.
Thus we have, for $x \in \calR$
\begin{eqnarray}  \label{preformula of sph}
\omega_\nu(x; s) = 
\frac{\gamma(\sigma(\chi))}{Q \cdot c_\sigma(\chi)} \cdot
\dsum{\sigma \in W_0} \calP_U(\wt{T}_\sigma^\chi(\Psi_\nu(x,s; \; )))(1).
\end{eqnarray}

\bigskip
On the other hand, by Lemma~\ref{lemma 1.8}, there exists an invertible matrix $A_\sigma(\chi)$ for $\sigma \in W_0$ satisfying
\begin{eqnarray}   \label{feq of T}
\left( \wt{T_\sigma^\chi}(\Psi_\nu(x, s; \; )) \right)_{\nu \in \calU} =
A_\sigma(\chi) \big{(} \Psi_\nu(x, \sigma(s) + \ve_\sigma; \; ) \big{)}_{\nu \in \calU},
\end{eqnarray}
where $A_\sigma(\chi)$ depends only on the $G$-orbit containing $x$, since $\wt{T_\sigma^\chi}$ is $G$-equivariant and (\ref{G-action on Psi}).

For $x \in \calR$, we obtain by (\ref{preformula of sph}) and (\ref{feq of T})
\begin{eqnarray*}
\left( \omega_\nu(x;s ) \right)_{\nu \in \calU} 
&=&
\frac{1}{Q} \dsum{\sigma \in W_0}\, \frac{\gamma(\sigma(\chi))}{c_\sigma(\chi)} \cdot 
A_\sigma(\chi) \big{(} \calP_U(\Psi_\nu(x,\sigma(s)+\ve_\sigma; \; )))(1) \big{)}_{\nu} \\  
&=&
\frac{1}{Q} \dsum{\sigma \in W_0}\, \gamma(\sigma(\chi)) \cdot B_\sigma(\chi) 
\left( \wt{\Psi}_\nu(x,\sigma(s)+\ve_\sigma; 1) \right)_{\nu},
\end{eqnarray*}
where we set $B_\sigma(\chi) = c_\sigma(\chi)^{-1} A_\sigma(\chi)$.
By (\ref{isom}) and Proposition~\ref{extension of T}, we see the invertible matrix $B_\sigma(\chi)$ satisfies the cocycle relation
$$
B_{\sigma\tau}(\chi) = B_{\tau}(\chi) B_\sigma(\tau(\chi)), \quad \sigma, \tau \in W_0.
$$
Hence
\begin{eqnarray}
\left( \omega_\nu(x;\tau(s)+\ve_\tau) \right)_\nu &=& 
\frac{1}{Q} \dsum{\sigma \in W_0} \gamma(\sigma\tau(\chi)) B_\sigma(\tau(\chi)) \left( \wt{\Psi}_\nu(x, \sigma\tau(s)+\ve_{\sigma\tau}; 1) \right)_\nu \nonumber \\
&=&
B_{\tau}(\chi)^{-1} \left( \omega_\nu(x;s) \right)_\nu,  \label{fun eq of sph}
\end{eqnarray}
and the above relation determines $B_\sigma(\chi)$, since $\set{\omega_\nu(x;\sigma(s)+\ve_\sigma)}{\nu \in \calU}$ is  linearly independent for each $\sigma \in W_0$. 

Finally, if $x \in \calR^+$, we
have $\wt{\Psi}_\nu(x, s; 1) = \abs{f(x)}^s$.
\qed

\bigskip
Now we assume the condition (A.4) for $(\B, \X)$ and recall $\wt{\omega}_u(x;s)$, 
then we do not need to consider the subgroup $W_0$, i.e. $W_0 = W$. 
We have $\chi = \chi_s = \abs{\psi}^s$ for $\wt{\omega}_u(x; s)$ (cf. (\ref{chi_s}), (\ref{def Psi_nu})), and instead of Theorem~\ref{Th1}, we have the following.

\begin{thm}
Assume (A1), (A2), (A3) and (A4) for $(\B, \X)$ and $s$ is generic.
For $x \in \calR$, one has 
$$
\big{(} \wt{\omega}_\nu(x; s) \big{)}_{\nu \in \calU} = 
\frac{1}{Q} \dsum{\sigma \in W}\, \gamma(\sigma(s)) \cdot B_\sigma(s) \cdot 
\left( \dint{U}\, \abs{f(u\cdot x)}_\nu^{\sigma(s)}du \right)_{\nu \in \calU}.
$$
Moreover, if $x \in \calR^+$, one has
$$
\big{(} \wt{\omega}_\nu(x; s) \big{)}_{\nu \in \calU} = 
\frac{1}{Q} \dsum{\sigma \in W}\, \gamma(\sigma(s))) \cdot B_\sigma(s) \cdot 
\left( \abs{f(x)}_\nu^{\sigma(s)} \right)_{\nu \in \calU}.
$$
Here the constant $Q$ and the rational function $\gamma(s) = \gamma(\abs{\psi}^s)$ of $q^{s_1},\ldots, q^{s_n}$ are determined by the group $G$ as in \S 2.1, and the matrix $\wt{B}_\sigma(s)$ is determined by
the functional equation
$$
\big{(} \wt{\omega}_\nu(x; s) \big{)}_{\nu \in \calU} = 
\wt{B}_\sigma(s) \big{(} \wt{\omega}_\nu(x; \sigma(s)) \big{)}_{\nu \in \calU}.
$$ 
\end{thm}

\vspace{2cm}
\Section{}  %\S 3
We follow the previous notations, and assume that $(\B, \X)$ satisfies (A1), (A2) and (A3).
In this section we give a condition to assure the existence of functional equations for $\sigma \in W$ attached to a simple root $\alp$, and explain how the functional equations are reduced to those of $p$-adic local zeta functions of small prehomogeneous vector spaces of limited type. A basic reference is \cite{Hamb}, where we assumed (A4) also.

\bigskip
\noindent
{\bf 3.1.} 
For a simple root $\alp$ whose associated reflection $\sigma = \sigma_\alp$ belongs to $W_0$, denote by $\bP$ the standard parabolic subgroup $\bP_{\{ \alp \} }$ in the sense of \cite{Bo}-21.11, and consider the following condition.

\medskip
\noindent
(A5) There exists a $k$-rational representation $\rho : \bP \longrightarrow R_{k'/k}(GL_2)$  satisfying  
\begin{eqnarray} \label{assum:rho}
&& 
\rho(\bP) = R_{k'/k}(GL_2) \mbox{ or }R_{k'/k}(SL_2), \quad 
\rho(\sigma) = \twomatrix{0}{1}{-1}{0} \big{(}= {\bf j}, \mbox{ say } \big{)},  \nonumber \\
&&
\rho^{-1}(\B_2) \subset \B, 
\quad \rho(K \cap \bP) \supset R_{k'/k}(SL_2)(\calO), 
\end{eqnarray}
where $k'$ is a finite unramified extension of $k$, $R_{k'/k}$ is the restriction functor of base field, and $\B_2$ is the Borel subgroup of $\rho(\bP)$ consisting of upper triangular matrices. 

\medskip
Chevalley groups are typical examples which have $\rho$ as above for $k = k'$ for each simple root (cf. \cite{Sf2}-\S 4.1.). As for $R_{k'/k}$, we note an example in \S 4.3.

\medskip
For each $u \in J(X)$, we set $J_u = \set{\nu \in J(X)}{P \cdot X_\nu = P \cdot X_u}$.
Denote by $e$ the group index $[\frX(\B) \cap (\frX_0(\B)\otimes_\Z \Q) : \frX_0(\B)]$, by $d$ the extension degree of $k'/k$, and let $\ve_\sigma$ be the same as in (\ref{ve-sigma}). 
Our first main result of this section is the following (cf. \cite{Hamb}-Theorem~2.6).

\begin{thm} \label{fun-eq for omega_u}
We assume (A1), (A2), (A3) for $(\B, \X)$. For a simple root whose associated reflection $\sigma$ belongs to $W_0$, we assume (A5), and keep the notations above. Then, there exists a functional equation 
\begin{eqnarray*}
\omega_u(x; s) = \frac{1 - q^{-2d-\sum_i e_i s_i}}{1 - q^{-2d-\sum_i e_i (\sigma(s)_i+\ve_i)}} \times
\sum_{\nu \in J_u}\, \gamma_{u\nu}(s) \cdot \omega_\nu(x; \sigma(s)+\ve_\sigma),
\end{eqnarray*}
where $\ve_i$ is the $i$-th component of $\ve_\sigma$, $\gamma_{u\nu}(s)$'s are rational functions of $q^{\frac{s_i}{e}}$, and $e_i$ is defined in (\ref{e_i}).
\end{thm}

\medskip
Hereafter we assume $(\B, \X)$ satisfies (A1), (A2) and (A3) and $\bP$ satisfies (A5).
In order to prove the above theorem and explain about gamma-factors, we introduce the following space.

Set $\wt{\X} = \X \times \V$ with $\V = R_{k'/k}(M_{21})$ and $\wt{\bP} = \bP \times R_{k'/k}(GL_1)$, and define the action  
\begin{eqnarray}
(p,t) \cdot (x,v) = (p\cdot x, \rho(p)vt^{-1}), \quad (p,t) \in \wt{\bP}, \; (x,v) \in \wt{X}.
\end{eqnarray}
Here we identify $k'$ with its image by the regular representation in $M_d(k)$ (with respect to a fixed basis for $k'/k$) and realize $R_{k'/k}(GL_2)$ (resp. $\V$) in $GL_{2d}(\overline{k})$ (resp. $M_{2d,d}(\overline{k})$), where $\overline{k}$ is the algebraic closure of $k$. Then we may identify as $\wt{P} = P \times GL_1(k')$ and $V = {k'}^2$. 
Further we regard $\B$ as a subgroup of $\wt{\bP}$ by the embedding
\begin{eqnarray}
\B \longrightarrow \wt{\bP}, \; b \longrightarrow (b, \rho(b)_1),
\end{eqnarray}
where $\rho(b)_1$ is the upper left $d$ by $d$ block of $\rho(b) \in R_{k'/k}(GL_2)$. Then one can identify $\B$ as the stabilizer subgroup of $\wt{\bP}$ at $v_0 = \twovector{1}{0} \in V$, i.e.
\begin{eqnarray}
\B \cong \wt{\bP}_{v_0} = \set{(p,t) \in \wt{\bP}}{\rho(p)v_0t^{-1} = v_0}.
\end{eqnarray}
 
Then we have the following(cf. \cite{Hamb}-Lemma~1.1, Proposition~1.2).

\begin{prop} \label{S3-Prop1}
{\rm (i)} One has the following isomorphism: 
$$
\begin{array}{rll}
\frX(\wt{\bP}) \cong \frX(\bP) \times \frX(R_{k'/k}(GL_1)) & \longrightarrow & \frX(\B)\\
(\psi_1, \psi_2) & \longmapsto & [p \mapsto \psi_1(p)\psi_2(\rho(p)_1)] .
\end{array}
$$

\noindent
{\rm (ii)} The space $(\wt{\bP}, \wt{\X})$ satisfies (A1), (A2) and (A3).  Further, if $(\B, \X)$ satisfies also (A4), then so does $(\wt{\bP}, \wt{\X})$, and $[\frX(\B) : \frX_0(\B)] = [\frX(\wt{\bP}) : \frX_0(\wt{\bP})]$.

\noindent
{\rm (iii)} The set of open $B$-orbits in $X$ corresponds bijectively to the set of open $\wt{P}$-orbits in $\wt{X}$ by the map $B\cdot x \mapsto \wt{P}\cdot(x,v_0)$.
\end{prop}

Let $\set{\wt{f}_i(x,v)}{1 \leq i \leq n}$ be the basic set of relative $\wt{\bP}$-invariants, which are regular on $\wt{\X}$ and satisfy $f_i(x) = \wt{f}_i(x, v_0)$. 
Since $\wt{f}_i(x,v)$ is a relative $R_{k'/k}(GL_1)$-invariant with respect to the action on $v$, it is homogeneous in the coordinates of $v$ over $k$, and we set
\begin{eqnarray} \label{e_i}
e_i = \deg_v \wt{f}_i(x,v),  \quad 1 \leq i \leq n.
\end{eqnarray}

We denote by $\wt{\psi}_i$ the character corresponding to $\wt{f}_i(x,v)$, then $\psi_i = \wt{\psi}_i \vert_\B$ for each $i$. 
For each $u \in J(X)$ denote by $\wt{X_u}$ the $\wt{P}$-orbit corresponding to $X_u$ by Proposition~\ref{S3-Prop1}, then we have
$$
X^{op} = \dsqcup{u \in J(X)}\, X_u, \qquad \wt{X}^{op} = \dsqcup{u \in J(X)}\, \wt{X}_u.
$$
Further we see
\begin{eqnarray}  \label{orbit decomp}
\wt{X}_u = \dsqcup{\ell \in \Z}\, (1, \pi^\ell) \cdot \wt{X_{u,0}}, \quad
\wt{X_{u,0}} = \dsqcup{h \in SL_2(\calO')/\Gamma}\, (\wt{h},1) \cdot(X_u \times \{v_0 \}),
\end{eqnarray}
where $\calO'$ is the ring of integers in $k'$, $\wt{h} \in K \cap \bP$ satisfying $\rho(\wt{h}) = h$ for each $h \in SL_2$ and $\Gamma = \set{\twomatrix{1}{a}{0}{1} \in SL_2(\calO')}{a \in \calO'}$.

Denote by $\SX$ and $\SXtilde$ the spaces of Schwartz-Bruhat functions on $X$ and $\wt{X}$, respectively. 
For $s \in \C^n$ and $u \in J(X)$, we consider the following integrals, which we call zeta integrals,
\begin{eqnarray*} 
&&
\Omega_u(\phi; s) 
=
\dint{X}\, \phi(x)\cdot \abs{f(x)}_u^{s} dx, \qquad \phi \in \SX, \\
&&
\wt{\Omega}_u(\wt{\phi}; s) 
= 
\dint{\wt{X}}\, \wt{\phi}(x,v) \cdot \abs{\wt{f}(x,v)}_u^{s} dxdv,
\qquad \wt{\phi} \in \SXtilde,
\end{eqnarray*}
where $dx$ is a $G$-invariant measure on $X$, $dv$ is a Haar measure on $V$, $\abs{f(x)}_u^s$ is defined in (\ref{def:sph}), and
\begin{eqnarray*} \label{def:tilde f^s}
\abs{\wt{f}(x,v)}^s = \ddprod{i=1}{n}\, \abs{\wt{f}_i(x,v)}^{s_i}, \quad
\abs{\wt{f}(x,v)}_u^s = 
\left\{ \begin{array}{ll}
\abs{\wt{f}(x,v)}^{s} & \mbox{if } (x,v) \in \wt{X_u}\\
0 & \mbox{ otherwise }. 
\end{array} \right.
\end{eqnarray*}
The above integrals are absolutely convergent for $\real(s_i) \geq 0, \; 1 \leq i \leq n$, and analytically continued to rational functions of $q^{s_i}, \; 1 \leq i \leq n$. We have the following,
where the assertion (i) is clear, and the assertion (ii) follows from (\ref{orbit decomp}).

\begin{lem} \label{relations}
{\rm (i)} Let $ch_x$ be the characteristic function of $K \cdot x$ in $\calS(X)$, then
\begin{eqnarray} \label{omega-Omega}
\omega_u(x; s) = v(K \cdot x)^{-1} \cdot \Omega_u(ch_x; s), \qquad x \in X, \; u \in J(X),
\end{eqnarray}
where $v(K\cdot x)$ is the volume of $K \cdot x$ by the above measure $dx$.
 
{\rm (ii)} Let $\wt{\phi} = \phi \otimes ch_{V(m)}$, where $\phi \in \SKX$ and $ch_{V(m)}$ is the characteristic function of $V(\pi^m\calO_{k'})$ in $\calS(V)$. Then 
\begin{eqnarray}
\wt{\Omega}_u(\wt{\phi}; s) = c \cdot \frac{q^{-m(2d+\sum_i e_is_i)}}{1 - q^{-2d-\sum_ie_is_i}} \times
\Omega_u(\phi; s), \quad u \in J(X),
\end{eqnarray}
where $c$ is a constant depending only on the normalization of measures, and independent of $u$.
\end{lem}

\bigskip
In order to study the action of $\sigma$ on $s$ for $\wt{\Omega}_u(\wt{\phi};s)$, we define the partial Fourier transform $\calF$ on $\calS(\wt{X})$ by
\begin{eqnarray*}
\calF(\wt{\phi})(x,v) = \dint{V}\, \eta({}^tv {\bf j} w) \wt{\phi}(x,w)dw, 
%\qquad ( \, {\bf j} = \twomatrix{0}{1}{-1}{0}),
\end{eqnarray*}
where $\eta$ is an additive character on $k'$ of conductor $\ell$, and consider the following distribution on $\calS(\wt{X})$
\begin{eqnarray*}
T_{u,s}(\wt{\phi}) = \wt{\Omega}_u(\wt{\phi}; s), \quad T^*_{u,s}(\wt{\phi}) = T_{u,s}(\calF(\wt{\phi})).
\end{eqnarray*}

We examine the relative invariancy of these distribution concerning the action of $\wt{P}$ on $\calS(\wt{X})$, where ${}^{\wt{p}} \wt{\phi}(x,v) = \wt{\phi}(\wt{p}^{-1} \cdot (x,v)), \; \wt{p} \in \wt{P}$.
We note the action of $\sigma$ on characters in the following (cf. \cite{Hamb}-Lemma 2.1).
\begin{lem} \label{lem char}
{\rm (i)} For a character $\xi \in \frX(\wt{\bP})$,
$$
\sigma(\xi)(p,t) = \xi(p, \frac{\det \rho(p)}{t}), \qquad (p,t) \in \wt{P}.
$$
{\rm (ii)} 
$
\abs{\wt{\psi}(p,t)}^{\ve_\sigma} = \abs{N_{k'/k}(t^2\det \rho(p)^{-1})}, \qquad (p,t) \in \wt{P}.
$
\end{lem}

\begin{prop}
The distributions $T_{u,s}^*$ and $T_{u, s^*}$ with $s^* = \sigma(s) + \ve_\sigma$ have the same relative invariancy with respect to the action of $\wt{P}$.
\end{prop}

\proof
First we obtain 
\begin{eqnarray*}
T_{u,s}({}^{(p,t)}\wt{\phi}) &=& \dint{X \times V}\, \wt{\phi}(p^{-1}\cdot x, \rho(p)^{-1}vt)\abs{\wt{f}(x,v)}^s dxdv\\
&=&
\abs{N_{k'/k}(t^{-2}\det \rho(p))} \dint{X \times V}\, \wt{\phi}(x,v) \abs{\wt{f}(p\cdot x, \rho(p)vt^{-1})}^s dx dv\\
&=&
\abs{\wt{\psi}(p,t)}^{s-\ve_\sigma} T_{u,s}(\wt{\phi}), 
\end{eqnarray*}
where we use Lemma~\ref{lem char}(ii) and $G$-invariancy of $dx$. Next, since ${\bf j} \rho(p) = \det\rho(p)\, {}^t\rho(p)^{-1} {\bf j}$, we have
\begin{eqnarray*}
\calF({}^{(p,t)}\wt{\phi})(x,v) &=& \dint{V}\, \eta({}^tv {\bf j} w) \wt{\phi}(p^{-1}\cdot x, \rho(p)^{-1}wt) dw\\
&=& 
\abs{\wt{\psi}(p,t)}^{-\ve_\sigma} \dint{V} \eta({}^t(\rho(p)^{-1}v {\scriptsize \frac{\det\rho(p)}{t}}) {\bf j}w) \wt{\phi}(p^{-1}\cdot x, w) dw\\
&=&
\abs{\wt{\psi}(p,t)}^{-\ve_\sigma} \left( {}^{(p, \frac{\det\rho(p)}{t})}\calF(\wt{\phi}) \right)(x,v).
\end{eqnarray*}
By the above calculation together with Lemma~\ref{lem char} (i), we obtain
\begin{eqnarray*}
T^*_{u,s}({}^{(p,t)}\wt{\phi}) = \abs{\wt{\psi}(p,t)}^{-\ve_\sigma + \sigma(s-\ve_\sigma)} T^*_{u,s}(\wt{\phi}) = \abs{\wt{\psi}(p,t)}^{\sigma(s)}T^*_{u,s}(\wt{\phi}).
\end{eqnarray*}
\qed
 
\bigskip
Because of the uniqueness of the relatively invariant distribution on homogeneous space(cf. \cite{Ig}-Proposition 7.2.1), we have the following identity 
\begin{eqnarray} \label{feq for T}
T_{u,s}^*(\wt{\phi}) = \dsum{\nu \in J(X)}\, \gamma^{\eta}_{u,\nu}(s) T_{\nu,s^*}(\wt{\phi}), \quad \wt{\phi} \in \calS(\wt{X}^{op}),
\end{eqnarray}
where $\gamma_{u,\nu}^\eta(s)$ is a constant independent of $\wt{\phi}$.  
Since $T_{u,s}(\wt{\phi})$ and $T^*_{u,\nu}(\wt{\phi})$ are continued to rational functions of $q^{s_1}, \ldots, q^{s_n}$ and $[\frX(\wt{\bP})\cap (\frX_0(\wt{\bP})\otimes_\Z \Q) : \frX_0(\wt{\bP})] = e$, the above $\gamma_{u,\nu}^\eta(s)$ are rational functions of $q^{s_1/e}, \ldots, q^{s_n/e}$. 

On the other hand,
under our assumption, essentially by (A1) and (A3), it is known that 
there is no nonzero distribution for generic $s$ whose support is contained in $\wt{X} \backslash \wt{X}^{op}$ and relative invariancy for $\wt{P}$ is $\abs{\wt{\psi}}^s$, where 'generic' means the same as in Lemma~\ref{lemma 1.8} (cf. \cite{Hamb}-(F5), \cite{Sf3}-Lemma~2.3, Corollary~2.4). 
Hence the identity (\ref{feq for T}) holds for any $\wt{\phi} \in \calS(\wt{X})$. Finally, if $\wt{\phi} \in \calS(\wt{X})$ is zero outside of $P \cdot X_u \times V$, then so does $\calF(\wt{\phi})$, hence we see 
$$
\gamma_{u,\nu}^\eta(s) = 0 \quad \mbox{ unless } P\cdot X_u = P\cdot X_\nu,\; \mbox{ i.e.,}\, \nu \in J_u.
$$
Thus we obtain the following theorem.

\begin{thm}  \label{th: tilde feq}
There exist rational functions $\gamma_{u\nu}^\eta(s)$ of $q^{\frac{s_1}{e}},\ldots, q^{\frac{s_n}{e}}$, which satisfy the following functional equation : 
\begin{eqnarray*} \label{feq-Omega}
\wt{\Omega}_u(\calF(\wt{\phi}); s) = \dsum{\nu \in J_u} \gamma_{u \nu}^\eta(s) \cdot
\wt{\Omega}_\nu(\wt{\phi}; \sigma(s)+\ve_\sigma),
\qquad \wt{\phi} \in \calS(\wt{X}).
\end{eqnarray*}
\end{thm}

\slit
We note here that $\gamma_{u\nu}^\eta(s)$ depends on the choice of the character $\eta$ and the normalization of $dv$ on $V$, since $\calF(\wt{\phi})$ does. Let normalize $dv$ on $V$ to be self dual with respect to the inner product $(v,w) \longmapsto \eta({}^tv{\bf j} w)$, so $vol(\V(\calO)) = q^{\ell d}$.

\begin{cor} \label{cor}
For any $\phi \in \SKX$, we have
\begin{eqnarray*}
\Omega_u(\phi;s) = 
\frac{ 1 - q^{-2d-\sum_{i}\, e_is_i} }
     { 1 - q^{-2d-\sum_{i}\, e_i(\sigma(s)_i + \ve_i)} } 
\times \dsum{\nu \in J_u}\, \gamma_{u \nu}(s) \cdot \Omega_\nu(\phi; \sigma(s)+\ve_\sigma),
\end{eqnarray*}
where $\ve_i$ is the $i$-th component of $\ve_\sigma$ and 
$$
\gamma_{u\nu}(s) = q^{\ell(d + \sum_i\, e_is_i)}\cdot \gamma_{u \nu}^\eta(s),
$$
which is independent of the choice of the character $\eta$ on $k'$.
\end{cor}

Now Theorem~\ref{fun-eq for omega_u} follows from Corollary~\ref{cor} and Lemma~\ref{relations}.

\bigskip
\noindent
{\bf 3.2.} 
In this subsection we look at $\V = R_{k'/k}(M_{21})$ together with the action of $(\rho(\bP_x) \times R_{k'/k}(GL_1))$ for $x \in X^{op}$. First observation is the following (\cite{Hamb}-Lemma~3.1).

\begin{lem} \label{lem:3-1}
For each $x \in X^{op}$, $(\rho(\bP_x) \times R_{k'/k}(GL_1), \V)$ is a prehomogeneous vector space defined over $k$. Further, for $v \in V$, 
$\rho(P_x)v{k'}^\times$ is open in $V$ if and only if $\wt{P}\cdot (x,v)$ is open in $\wt{X}$.
\end{lem}

For each $u \in J(X)$, fix an element $x_u \in X_u$ and denote by $\bP_u$ the stabilizer of $x_u$ in $\bP$. Then we obtain (\cite{Hamb}-Lemma~3.2)

\begin{lem} \label{lem:3-2}
{\rm (i)} 
For any $u, \nu \in J(X)$, prehomogeneous vector spaces $(\bP_\nu \times R_{k'/k}(GL_1), \V)$ $(\bP_u \times R_{k'/k}(GL_1), \V)$ are isomorphic. If $\nu \in J_u$, they are isomorphic over $k$.

\noindent
{\rm (ii)} 
The set of $k$-rational points of the open orbit in $(\rho(\B_u) \times R_{k'/k}(GL_1), \V)$ decomposes as
\begin{eqnarray} \label{k-orbit decomp}
\left(  \rho(\bP_u)v_0R_{k'/k}(GL_1) \right) (k) = \dsqcup{\nu \in J_u}\, \rho(P_up_\nu)v_0{k'}^\times,
\end{eqnarray}
where $p_\nu \in P$ satisfying $p_\nu^{-1} \cdot x_u \in X_\nu$.
\end{lem}

For $\wt{\phi} = \phi_1 \otimes \phi_2$ with $\phi_1 \in \calS(X)$ and $\phi_2 \in \calS(V)$, we have
\begin{eqnarray*}
\calF(\wt{\phi}) = \phi_1 \otimes \calF_V(\phi_2), 
\end{eqnarray*}
where
\begin{eqnarray*}
\calF_V(\phi_2)(v) = \dint{V}\, \eta({}^tv{\bf j} w)\phi_2(w) dw.
\end{eqnarray*}
Taking these $\wt{\phi}$ in Theorem~\ref{th: tilde feq} and pulling out the $V$-part, we obtain the following theorem (cf. \cite{Hamb}-Theorem~3.3) which shows that the functional equations of spherical functions $\omega_u(x; s)$ are reduced to those for "small" prehomogeneous vector spaces, and oppositely gamma factors can be calculated from those of these prehomogeneous vector spaces.

\begin{thm} \label{thm 3-1}
The prehomogeneous vector space $(\bP_u\times GL_1, \V)$ has the following functional equation: 
\begin{eqnarray*} \label{V-fun eq 2}
\lefteqn{\dint{V}\, \calF_V(\phi)(v)\abs{\wt{f}(x_u,v)}_u^{s} dv}\\ 
&=& 
\dsum{\nu \in J_u}\, \gamma_{u \nu}^\eta(s) \dint{V}\, \phi(v)\abs{\wt{f}(x_u,v)}_\nu^{\sigma(s)+\ve_\sigma} dv, 
\qquad \phi \in \calS(V), \nonumber
\end{eqnarray*}
where the gamma factors $\gamma_{u\nu}^\eta(s)$ are the same as those for $\wt{\Omega}_u(\wt{\phi}; s^*)$ in Theorem~\ref{th: tilde feq}.
\end{thm}

\mslit
Because of the existence of the functional equations of the above type, we see the following (cf. \cite{Hamb}-Theorem~3.6).

\begin{thm} \label{thm 3-2}
For the prehomogeneous vector space $(\rho(\bP_u) \times R_{k'/k}(GL_1), \V)$,
the identity component of $\rho(\bP_u) \times R_{k'/k}(GL_1)$ is isomorphic to $R_{k'/k}(GL_1 \times GL_1)$ over the algebraic closure $\overline{k}$ of $k$.
\end{thm}

\bigskip
\noindent
{\bf 3.3.} 
If the condition (A4) is also satisfied by $(\B, \X)$, we should consider $\wt{\omega}_u(X;s)$.
In this subsection, we assume (A1),(A2), (A3) and (A4) for $(\B, \X)$. We assume (A5) for a simple root $\alp$ whose associated reflection $\sigma$, and keep the notations before.
Then we have Theorem~\ref{mdf th: tilde feq} instead of Theorem~\ref{th: tilde feq}, and based on it, we obtain Theorem~\ref{mdf fun-eq for omega_u} and Theorem~\ref{mdf thm 3-1}, these are the original formulation in \cite{Hamb}. We do not need to modify Theorem~\ref{thm 3-2}. 

\begin{thm} \label{mdf th: tilde feq}
There exist rational functions $\wt{\gamma}_{u\nu}^\eta(s)$ of $q^{\frac{s_1}{e}},\ldots, q^{\frac{s_n}{e}}$, which satisfy the following functional equation : 
\begin{eqnarray*} \label{mdf feq-Omega}
\wt{\Omega'}_u(\calF(\wt{\phi}); s) = \dsum{\nu \in J_u} \wt{\gamma}_{u \nu}^\eta(s) \cdot
\wt{\Omega'}_\nu(\wt{\phi}; \sigma(s)),
\qquad \wt{\phi} \in \calS(\wt{X}),
\end{eqnarray*}
where
$$
\wt{\Omega'}_u(\wt{\phi};s) = \dint{\wt{X}}\wt{\phi}(x,v) \cdot \abs{\wt{f}(x,v)}_u^{s+\ve_0}dxdv.
$$ 
\end{thm}

\begin{thm} \label{mdf fun-eq for omega_u} 
Then exists a functional equation 
\begin{eqnarray*}
\wt{\omega}_u(x; s) = \frac{1 - q^{-2d-\sum_i e_i (s_i+\ve_i)}}{1 - q^{-2d-\sum_i e_i (\sigma(s)_i+\ve_i)}} \times
\sum_{\nu \in J_u}\, \wt{\gamma}_{u\nu}(s) \cdot \wt{\omega}_\nu(x; \sigma(s)),
\end{eqnarray*}
where $\ve_i$ is the $i$-th component of $\ve_0$, $\wt{\gamma}_{u\nu}(s)$'s are rational functions of $q^{\frac{s_i}{e}}$, and $e_i$ is defined in (\ref{e_i}).
\end{thm}

\begin{thm} \label{mdf thm 3-1}
The prehomogeneous vector space $(\bP_u\times GL_1, \V)$ has the following functional equation: 
\begin{eqnarray*} \label{V-fun eq 2}
\lefteqn{\dint{V}\, \calF_V(\phi)(v)\abs{\wt{f}(x_u,v)}_u^{s+\ve_0} dv}\\ 
&=& 
\dsum{\nu \in J_u}\, \wt{\gamma}_{u \nu}^\eta(s) \dint{V}\, \phi(v)\abs{\wt{f}(x_u,v)}_\nu^{\sigma(s)+\ve_0} dv, 
\qquad \phi \in \calS(V), \nonumber
\end{eqnarray*}
where the gamma factors $\wt{\gamma}_{u\nu}^\eta(s)$ are the same as those for $\wt{\Omega'}_u(\wt{\phi}; s)$ in Theorem~\ref{mdf th: tilde feq}.
\end{thm}

\begin{rem} \label{thm for sph with char}
{\rm 
We recall Remark~\ref{rem sph with char}. Assume $\calU$ is a subset containing $J_u$ and canonically
identified  with a subgroup of $k^{\times n}/\prod_i \psi_i(B)$, and denote by $\what{\calU}$ the character group of $\calU$. We may define similarly $\Omega(\phi; \chi; s)$ and $\wt{\Omega}(\wt{\phi};\chi;s)$ for $\chi \in \what{\calU}$. Then, instead of Theorem~\ref{fun-eq for omega_u}, we have
\begin{eqnarray}
\omega(x; \chi; s) = 
\frac{ 1 - q^{-2d-\sum_{i}\, e_i(s_i)}  }
     { 1 - q^{-2d-\sum_{i}\, e_i (\sigma(s)_i+\ve_i)} } 
\times \dsum{\xi \in \what{\calU}}\, A_{\chi \xi}(s) \omega(x; \xi; \sigma(s)+\ve_\sigma) ,
\end{eqnarray}
where
\begin{eqnarray*}
&&
A_{\chi \xi}(s) = \frac{1}{\sharp(\calU)} \dsum{u, \nu \in \calU}\, \chi(u) \overline{\xi}(\nu) \gamma_{u\nu}(s), \quad \gamma_{u\nu}(s) = 0 \mbox{ unless } \nu \in J_u,\\
&&
\gamma_{u\nu}(s) = \frac{1}{\sharp(\calU)^2} \dsum{\chi, \xi \in \what{\calU}}\, 
\overline{\chi}(u) \xi(\nu) A_{\chi \xi}(s).
\end{eqnarray*}

We have a similar formula for $\wt{\omega}(x; \chi; s)$.
} % end \rm
\end{rem}

\vspace{2cm}

\Section{}  %\S 4   Examples
We prepare some notations. 
For a matrix $x \in M_n$, we denote by $d_i(x)$ is the determinant of upper left $i$ by $i$ block of $x$, by $x_{ij}$ the $(i,j)$-component of $x$, and $x_i = x_{ii}$.
We set $J_n = \left( \begin{array}{ccc}
0  & {} & 1\\
 {} & \gyakuddots & {}\\
 1 & {} & 0 \end{array} \right) \in GL_n$, i.e., the matrix whose anti-diagonal components are $1$ and $0$ elsewhere.
Set $\Lam_n = \set{\lam \in \Z^n}{\lam_1 \geq \cdots \geq \lam_n}$ and $\Lam_n^+ = \set{\lam \in \Lam_n}{\lam_n \geq 0}$.

In the following, we may take $K = \G(\calO)$, where $\calO$ is the ring of integers in $k$, each $\X$ is a symmetric space except the one in \S 4.5, which is a spherical homogeneous space, and $(\B, \X)$ satisfies (A1) and (A2). In each case, the open orbit $\X^{op}$ is given as the non-vanishing set of basic relative invariants, and one can consider spherical functions with character, since $J(X)$ has a canonical group structure (cf. Remark~1.4). 

\medskip
\noindent
{\bf 4.1.} The space of symmetric forms.  \\
$\G = GL_n, \X = \set{x \in \G}{{}^tx = x}, g \cdot x = gx{}^tg$.\\
$\B$ is the Borel group consisting of lower triangular matrices in $\G$. \\
$f_i(x) = d_i(x), \; \psi_i(p) = (p_1\cdots p_i)^2, 1 \leq i \leq n$, $W = W_0 \cong S_n$.\\
$J(X) \cong \left( k^\times/k^{\times2} \right)^n \left( \cong \left( \Z/4\Z \right)^n  \mbox{  if  } 2 \notin (\pi) \right)$. \\
$(\B, \X)$ satisfies (A1)--(A4).  

We may take representatives of $K \backslash X$ in $\calR^+$, since $Diag(a_1, \ldots, a_n) \in X$ with $v_\pi(a_1) \leq \cdots \leq v_\pi(a_n)$ is contained in $\calR^+$. 

The condition (A5) is satisfied for each simple root, indeed for the transposition $(\alp \; \alp+1), \; 1 \leq \alp \leq n-1$, we have
\begin{eqnarray*}
&&
\bP = \bP_\alp = \set{(p) \in \G}{p_{ij} = 0 \; \mbox{ unless } i \geq j \mbox{ or } (i,j) = (\alp, \alp+1)},\\
&&
\rho: \bP \longrightarrow GL_2, \; p \longmapsto \twomatrix{p_{\alp+1,\alp+1}}{-p_{\alp+1,\alp}}{-p_{\alp,\alp+1}}{p_{\alp,\alp}}.
\end{eqnarray*}
The small prehomogeneous vector spaces are of type $(O(T) \times GL_1, V)$ for some symmetric matrix $T$ of size 2. 
 
Functional equations with respect to the Weyl group have been known by a different method based on the explicit expressions of spherical functions of size 2 (\cite{H1}-III), and one can apply Theorem~\ref{Th1} formally, but good expressions of spherical functions are not known for general $n \geq 3$, only partial results are known(\cite{H1}).

\bigskip
\noindent
{\bf 4.2.} The space of alternating forms.\\
$\G = GL_{2n}, \X = \set{x \in \G}{{}^tx = -x}, g \cdot x = gx{}^tg$.\\
$\B = $ the Borel group consisting of lower triangular matrices in $\G$.\\
$f_i(x) = {\rm pf}_i(x), \; \psi_i(p) = p_1\cdots p_{2i}, \; 1 \leq i \leq n$, where ${\rm pf}_i(x)$ is the phaffian of the upper left $2i$ by $2i$ block of $x$.\\
$W = S_{2n} \supset W_0 \cong S_n$, in fact
for each $\sigma \in S_n$, we associate $w_\sigma \in W$ such that $w_\sigma(2i-1) = 2\sigma(i)-1, \; w_\sigma(2i) = 2\sigma(i), \; 1 \leq i \leq n$.\\
$X^{op}$ is a single $B$-orbit.
As a set of complete representatives of $K \backslash X$, we may take 
$$
\set{\pi^\lam}{\lam \in \Lam_n}, \quad \pi^\lam = \twomatrix{0}{\pi^{\lam_1}}{-\pi^{\lam_1}}{0} \bot \cdots \bot \twomatrix{0}{\pi^{\lam_n}}{-\pi^{\lam_n}}{0} \in X,
$$
and $J_{2n} \cdot \pi^\lam \in \calR^+$.
 
The explicit formula of $\omega(x; s)$ was calculated by another method (\cite{HS1}), and it can be reproduced also by Theorem~2.5. We introduce a new variable $z$ related to $s$ by
$$
s_i = -z_i + z_{i+1} -2  \; \; (1 \leq i \leq n-1), \quad s_n = -z_n + n-1
$$
and write $\omega(x; z) = \omega(x; s)$. For each $\lam \in \Lam_n$
$$
\omega(\pi^\lam; s) = c_\lam \cdot \prod_{1 \leq  i < j \leq n }\, \frac{ 1 - q^{z_i-z_j-1}}{1 - q^{z_i-z_j+1} } \cdot P_\lam(q^{z_1}, \cdots, q^{z_n}; q^{-2}),
$$
where $c_\lam$ is an explicitly given constant in $\Q(q^{-1})$ and $P_\lam$ is a Hall-Littlewood symmetric polynomial (a symmetric Laurent polynomial of $q^{z_1},\ldots,q^{z_n}$). 
The Hall-Littllewood polynomial $P_\lam(x; t)$ is defined as follows (cf. \cite{Ma2})
\begin{eqnarray*}
&&
P_\lam(x; t) = P_\lam(x_1, \ldots, x_n; t) = \frac{(1-t)^n}{w_\lam(t)} \cdot \sum_{\sigma \in S_n}\, x_{\sigma(1)}^{\lam_1} \cdots x_{\sigma(n)}^{\lam_n} \prod_{i < j}\, \frac{x_{\sigma(i)}-tx_{\sigma(j)}}{x_{\sigma(i)}-x_{\sigma(j)}},\\
%\end{eqnarray*}
%where
%\begin{eqnarray*}
&&
w_\lam(t) = \prod_{j=1}^{r} \prod_{i=1}^{n_i} ( 1 - t^i), \quad \mbox{for } \lam = (\ell_1^{n_1}\cdots \ell_u^{n_r}), \; \ell_1 > \cdots > \ell_r, \; n_1 + \cdots + n_r = n,
\end{eqnarray*}
where the set $\set{P_\lam(x; t)}{\lam \in \Lam_n^+}$ forms a $\Z[t]$-basis for $\Z[t][x_1, \cdots, x_n]^{S_n}$, and the set\\
$\set{P_\lam(x; t)}{\lam \in \Lam_n}$ forms a $\Z[t]$-basis for $\Z[t][x_1^{\pm 1}, \cdots, x_n^{\pm 1}]^{S_n}$.

Setting 
$$
\Psi(x; z) = \omega(x; z) / \omega(\pi^{\bf 0} ; z), \qquad {\bf 0} = (0,\ldots,0) \in \Lam_n,
$$
we have the spherical transform which is a surjective $\hec$-module homomorphism
$$
\calS(K\backslash X) \longrightarrow \C[q^{\pm z_1}, \ldots, q^{\pm z_n}]^{S_n}, \; 
\phi \longmapsto \dint{X} \phi(x) \Psi(x; z) dx
$$
where $dx$ is $G$-invariant measure on $X$ and $\hec$ acts on the right hand side via $\lam_z(\phi) = \lam_s(\phi)$, a specialization of Satake transform $\hec \cong \C[q^{\pm t_1}, \ldots, q^{\pm t_{2n}}]^{S_{2n}}$.

Each spherical function on $X$ is associated to some $z \in \C^n/S_n$ through $\lam_z$, and it is a constant multiple of $\Psi(x; z)$.

\bigskip
\noindent
{\bf 4.3.} The space of hermitian forms.\\ 
For a quadratic extension $k'/k$ with involution $*$, we consider spherical functions on the space of hermitian forms $X = \set{x \in GL_n(k')}{x^* = x}$ with canonical action of $G = GL_n(k')$, where $(i,j)$-component of $g^*$ for $g = (g_{ij}) \in G$ is $g_{ji}^*$. We have to realize these objects as the sets of $k$-rational points: taking $u \in k'$ such $k'= k(u)$ and $u^2 \in k$, we identify $k'$ with the image of the inclusion
$$
k' \longrightarrow M_2(k), \; a+bu \longmapsto \twomatrix{a}{bu^2}{b}{a}
$$
and realize $\G = R_{k'/k}(GL_n)$ and the space $\X$ of hermitian forms in the following.
\begin{eqnarray*}
&&
\G = \set{(g_{ij}) \in GL_{2n}(\overline{k})}{g_{ij} = \twomatrix{a_{ij}}{b_{ij}u^2}{b_{ij}}{a_{ij}} \in M_2(\overline{k}) \quad (1 \leq i, j \leq n)},\\
&&
\X = \set{(x_{ij}) \in \G}{x_{ij} = \twomatrix{1}{0}{0}{-1} x_{ji} \twomatrix{1}{0}{0}{-1} \quad (1 \leq i, j \leq n)},\\
&&
\B = \set{g = (g_{ij}) \in \G}{g_{ij} = 0 \; \mbox{  unless } n \geq i \geq j \geq 1}.
\end{eqnarray*}
Then $G = \G(k) = GL_n(k'), X = \X(k) = \set{x \in G}{x^* = x}, g\cdot x = gxg^*$.\\
$B = \B(k)$ is the Borel subgroup of $G$ consisting of lower triangular matrices.\\
$f_i(x) = d_i(x) \in k$ for $x \in X$, and 
$\psi_i(p) = N_{k'/k}(p_1\cdots p_i)$ for $p \in B$, $1 \leq i \leq n$. \\
$W = W_0 \cong S_n$.  
$J(X) \cong \left( k^\times/N_{k'/k}(k^{' \times})  \right)^n \cong (\Z/2\Z)^n$.\\
$(\B, \X)$ satisfies (A1) -- (A4).

Functional equations with respect to the Weyl group have been known by a different method based on the explicit expressions of spherical functions of size 2 (\cite{H1}-III), and one applies Theorem~\ref{Th1} formally.

\medskip
Theorem~\ref{Th1} (of its original formulation) was used to obtain the explicit formula for the case unramified hermitian forms (\cite{JMSJ}-\S2). 
We consider the spherical function of type (\ref{char-mod-sph}) 
$$
\omega(x; z) = \omega(x; s) = \dint{K} \abs{f(k \cdot x)}^{s+\ve}dk,
$$
where $\ve = (-1,\ldots,-1,-\frac{n-1}{2}) + (\frac{\pi\sqrt{-1}}{\log q}, \ldots, \frac{\pi\sqrt{-1}}{\log q}) \in \C^n$, and $z$ is the new variable related to $s$ by 
$$
s_i = -z_i +z_{i+1} \quad (1 \leq  i\leq n-1), \quad s_n = -z_n.
$$
Though we are shifting the variable $s$ here, $\omega(x;z)$ is the same as before.
We have the functional equations 
\begin{eqnarray} \label{feq-uh}
\omega(x; z) = 
\prod_{\scriptsize 
\begin{array}{c}
1 \leq i < j \leq n\\
\sigma(i) > \sigma(j) \end{array} }\, 
\frac{q^{z_{\sigma(i)}} - q^{z_{\sigma(j)}-1}}{q^{z_{\sigma(j)}} - q^{z_{\sigma(i)}-1}} 
\times \omega(x; \sigma(z)), \quad \sigma \in S_n.
\end{eqnarray}
A set of complete representatives of $K \backslash X$ can be taken as
$$
\set{\pi^\lam = Diag(\pi^{\lam_1}, \ldots, \pi^{\lam_n})}{ \lam \in \Lam_n}
$$
and $J_n \cdot \pi^\lam \in \calR^+$. The explicit formula is given as 
$$
\omega(\pi^\lam; z) = c_\lam \cdot \prod_{1 \leq i < j \leq n}\, \frac{1 - q^{z_i- z_j- 1}}{1 + q^{z_i-z_j}} \cdot P_\lam(q^{z_1}, \ldots, q^{z_n};-q^{-1}),
$$
where $c_\lam$ is an explicitly given constant in $\Q(q^{-1})$ and $P_\lam$ is a Hall-Little symmetric polynomial (cf. \S 4.2). 
Setting
$$
\Psi(x; z) = \omega(x; z)/\omega(1_n; z),
$$
we have the spherical transform which is an $\hec$-module isomorphism
$$
\calS(K \backslash X) \cong \C[q^{\pm z_1}, \ldots, q^{\pm z_n}]^{S_n}, 
\phi \longmapsto \dint{X} \phi(x) \Psi(x; z) dx, 
$$
where $dx$ is the $G$-invariant measure on $X$ and $\hec$ acts on the right hand side via Satake transform $\wt{\lam}_z = \wt{\lam}_s : \hec \cong \C[q^{\pm 2z_1}, \ldots, q^{\pm 2z_n}]^{S_n}$. In particular $\calS(K \backslash X)$ is a free $\hec$-module of rank $2^n$. 
%(cf. \ref{modified sph}).   % Rem 1.3 

Each spherical function on $X$ is associated to some $z \in \C^n/S_n$ through $\wt{\lam}_z$, and the space of spherical functions associated to $\wt{\lam}_z$ has dimension $2^n$ and a basis\\ $\set{\Psi(x; z + \eps)}{\eps \in \{0, \frac{\pi\sqrt{-1}}{\log q} \}^n }$.

\medskip
For a simple root associated to the transposition $(\alp \; \alp+1), \; 1 \leq \alp \leq n-1$, there is a representation $\rho$ satisfying (A5), similarly given to the case of symmetric form, but we have to use $R_{k'/k}(GL_2)$ in order to define $\rho$ over $k$ (cf. \cite{Hamb}-\S 4.2). 
The small prehomogeneous vector spaces are of type $(U(T) \times GL_1, V)$ for some hermitian matrix $T$ of size 2. 

\medskip
Assume that $k'/k$ is ramified and take a prime element $\pi'$ of $k'$ as ${\pi'}^2 = -\pi$. Then some representatives of $K \backslash X$ contain a matrix of type $\twomatrix{0}{{\pi'}^{2m+1}}{-{\pi'}^{2m+1}}{0}$ as a direct summand, and they are not in $\calR^+$, functional equations of spherical functions are much more complicated (cf. \cite{H1}-III), and no good expressions of spherical functions are known for general $n \geq 3$.

\bigskip
\noindent
{\bf 4.4.} %Symmetric spaces of type $G(n) \big{/} G(r) \times G(n-r)$ with $2r \leq n$.\\
%$G(m) = GL_m, Sp_m, SO(m), U(m)$. $W = W(n)$ and $W_0 \cong W(r)$, where $W(m)$ is the Weyl group of $G(m)$ with respect to a maximal $k$-split torus for each group.

%\medskip
\noindent
(i) For a nondegenerate symmetric matrix $A$ of size $n$, we set 
$O(A) = \set{g \in GL_n}{A[g] = A}$ and $SO(A) = O(n) \cap SL_n$ where $A[g] = {}^tgAg$.
Set 
$$
\G = SO(H_n), \; H_n = \frac12 \twomatrix{}{1_n}{1_n}{}, \;  
\B = \set{\twomatrix{b_1}{b_2}{b_3}{b_4} \in \G}{\begin{array}{l}
           b_i \in M_n, \; b_3 = 0,\\
           b_1 \mbox{ is upper triangular} \end{array}}.
$$
In the following we consider the set of $k$-rational points. Fix a symmetric matrix $T \in GL_r(k)$, set\\
\hspace*{2cm}$\frX_T = \set{x \in M_{2n,r}}{H_n[x] = T}, \quad X_T = \frX_T/ O(T) \ni \overline{x} = xO(T)$,\\
and consider them as $\G$-spaces by left multiplication. The stabilizer of 
$$
\overline{x_T} \in X_T, \quad x_T = {}^t (T \; 0_{n-r,r} \; 1_r \; 0_{n-r,r}) \in \frX_T
$$
is isomorphic to $S(O(T) \times O(T \bot H_{n-r}))$, so the space $X_T$ is isomorphic to $SO(2n)/S(O(r) \times O(2n-r))$ over the algebraic closure of $k$.\\
A set of basic $B$-relative invariants and associated characters are given as follows
$$
f_{T, i}(\overline{x}) = d_i(T^{-1}[{}^tx_2]) = d_i(x_2T^{-1}{}^tx_2), \quad \psi_i(p) = (p_1\cdots p_i)^{-2}, \quad 1 \leq i \leq r,
$$
where $x_2$ is the submatrix of $x \in \frX_T$ consisting of its $(n+i)$-th row, $1 \leq i \leq r$, in order, and $p_i$ is the $i$-th diagonal component of $p \in B$. 
So $rank(\frX_0(B)) = r$, whereas $rank(\frX(B)) = n$.\\
$W = S_n \ltimes C_2^{n-1} \supset W_0 \cong S_r \ltimes C_2^{r-1}$, and $J(X_T) \cong \left(k^\times / k^{\times 2} \right)^{r-1}$. 

Functional equations with respect to the $S_r$-part are reduced to the case of symmetric forms.
When $r = n$, the condition (A5) is satisfied for each simple root and small prehomogeneous vector spaces are isomorphic to $(O(t)\times GL_1,V)$ for some symmetric matrix $t$ of size $2$, for details see \cite{Hamb}-\S4-3. The spherical functions on this space with respect to the Siegel parabolic subgroup have a close relation to Siegel singular series (cf. \cite{SH}).

\medskip
For odd $2n+1$, we may start with 
$$
\G = SO(H_n), \quad H_n = \left( \begin{array}{ccc}
                                         {} & {} & 1_n\\    
                                         {} & 1 & {}\\   
                                         1_n & {} & {}
                                  \end{array} \right) \in GL_{2n+1},
$$
and consider the space for symmetric $T \in GL_r(k)$                                  
$$
X_T = \frX_T/O(T), \quad \frX_T = \set{x \in M_{2n+1, r}}{H_n[x] = T}.
$$
Relative invariants are given similarly to the even case, and we have $W = S_n \ltimes C_2^n \supset W_0 \cong S_r \ltimes C_2^r$.

\bigskip
\noindent
(ii) For a nondegenerate hermitian matrix $A$ of size $n$, we set $U(A) = \set{g \in GL_n}{A[g] = g}$, where $A[g] = g^*Ag$ and $g^*$ is the same as in \S 4.3. 
In the following we write by the set of $k$-rational points for simplicity of notations (cf. \S 4.3).  Setting
\begin{eqnarray*}
&&
G = U(H_n), \; H_n = \twomatrix{}{1_n}{1_n}{}, B = \set{\twomatrix{b_1}{b_2}{b_3}{b_4} \in G}
{\begin{array}{l}
b_i \in M_n, \;; b_3 = 0\\
b_1 \mbox{ is upper triangular} \end{array} },
\end{eqnarray*}
define the $G$-spaces $\frX_T$ and $X_T = \frX_T/U(T)$ similarly to the case (i) for hermitian $T \in GL_r(k')$. Then, \\
$F_{T,i}(\overline{x}) = d_i(T^{-1}[x_2^*])$, $\psi_i(p) = N_{k'/k}(p_1\cdots p_i)^{-1}\;, 1 \leq i \leq r$, 
$J(X_T) \cong \left( k^\times/N_{k'/k}(k'^\times) \right)^{r-1}$.\\
$rank(\frX(\B)) = n, \; rank(\frX_0(\B)) = r$, and $W = S_n \ltimes C_2^n \supset W_0 = S_r \ltimes C_2^r$.\\
Functional equations with respect to the $S_r$-part are reduced to the case of hermitian forms, so we know well if $k'/k$ is unramified (cf. \S 4.3).

\medskip
Let us assume $r = n$ and $k'/k$ is unramified, and consider the spherical function on $X_T$ of type (\ref{char-mod-sph})
$$
\omega_T(\overline{x}; z) = \omega_T(\overline{x}; s) = \dint{K}\abs{f_T(kx)}^{s+\ve} dk, 
$$
where $\ve = (-1,\ldots,-1,-\frac{1}{2}) + (\frac{\pi\sqrt{-1}}{\log q}, \ldots, \frac{\pi\sqrt{-1}}{\log q}) \in \C^n$, and $z$ is the new variable related to $s$ by 
$$
s_i = -z_i +z_{i+1} \quad (1 \leq  i\leq n-1), \quad s_n = -z_n.
$$
This $\omega_T(\overline{x}; z)$ satisfies the same functional equation as in (\ref{feq-uh}) with respect to $S_n$, independent of the choice of $T$. Further we obtain 
\begin{eqnarray} \label{feq-unitary}
\omega_T(\overline{x}; z) = \abs{2}_k^{2z_n} \omega_T(\overline{x}; \tau(z)), \quad
\tau(z) = (z_1, \ldots, z_{n-1}, -z_n),
\end{eqnarray}
and we have functional equations with respect to whole $W$ by cocycle relations. The parabolic subgroup attached to $\tau$ does not have the representation satisfying (A5) and the above functional equation (\ref{feq-unitary}) does not come from prehomogeneous vector space.
Explicit formulas of $\omega_T(\overline{x})$ for some particular points are obtained by using these functional equations and Theorem~2.6.

Similar to the case (i), spherical functions on this space with respect to the Siegel parabolic subgroup have a close relation to hermitian Siegel series, for details see \cite{prep}.
 
%\bigskip
%\noindent
%(iii) For $G(n) = GL_n$ or $Sp_n$, we may consider the space isomorphic to $G(n) / G(r) \times G(n-r)%$ iver $k$ with $2r \leq n$. Then 
%$W = W(n)$ (the Weyl group of $G(n)$) and $W_0 \cong W(r)$. The open Borel orbit is a single orbit ev%enover $k$, and explicit formula is known by S.~Kato (\cite{Kato}) and Z.~Mao and S.~Rallis (unpublis%hed preprint).

%                                  
\bigskip
\noindent
{\bf 4.5.} $Sp_2 \times (Sp_1)^2$-space $Sp_2$ (cf. \cite{JNT}). \\
We assume that $k$ has odd residual characteristic. 

Setting
$$
Sp_n = \set{x \in GL_{2n}}{H_n[x] = H_n}, \; H_n = \twomatrix{}{1_n}{-1_n}{}, \; A[x] = {}^tx Ax,
$$
we embed $(Sp_1)^2 = (SL_2)^2$ into $Sp_2$ by 
$$
(\twomatrix{a}{b}{c}{d}, \twomatrix{e}{f}{g}{h}) \longmapsto 
\left( \begin{array}{cc|cc}   a & {} & b & {}\\
{} & e &  {} & f\\ 
\hline c & {} & d \\ {} & g & {} & h \\ \end{array}\right),
$$
and define the action of $\G = Sp_2 \times (Sp_1)^2$ on $\X = Sp_2$
$$
\wt{g} \cdot x = g_1x{}^tg_2, \quad \wt{g} = (g_1, g_2) \in \G, \; x \in \X.
$$
Then this space $\X$ is not a symmetric space, but a spherical homogeneous space, whereas 
$Sp_{2n}\big{/} (Sp_n \times Sp_n)$-space $Sp_{2n}$ is no longer spherical for $n \geq2$
 (and hence there is no open Borel orbit in it). 

We take the Borel subgroup of $\G$ consisting of matrices of type
$$
b = (
\left( 
 \begin{array}{c|c} 
\twomatrixminus{*}{*}{0}{*} & *\\
\hline
0 & \twomatrixminus{b_1}{0}{*}{b_2} \end{array}
\right),
\left( \begin{array}{c|c} 
\twomatrixminus{b_3}{0}{0}{b_4} & 0 \\
\hline
\twomatrixminus{*}{0}{0}{*} & \twomatrixminus{*}{0}{0}{*} \end{array}
\right)
 ) \in \G,
$$ 
then a set of basic relative $\B$-invariants on $\X$ and associated characters of $\B$ are given as 
\begin{eqnarray*}
\begin{array}{ll}
f_1(x) = x_{31}, & \psi_1(b) = b_1b_3,\\
f_2(x) = x_{32}, & \psi_2(b) = b_1b_4,\\
f_3(x) = x_{31}x_{42} - x_{32}x_{41}, & \psi_3(b) = b_1b_2b_3b_4,\\
f_4(x) = x_{31}x_{43}-x_{41}x_{33}, & \psi_4(b) = b_1b_2,\\
\end{array}
\end{eqnarray*}
where $x_{ij}$ is the $(i,j)$-component of $x$ and $b_i$ comes from the above expression of $b \in \B$.
$J(X) \cong k^\times/k^{\times 2}$, and $(\B, \X)$ satisfies (A1) -- (A4).

We consider the spherical function of type (\ref{char-mod-sph})
$$
\wt{\omega}(x; \chi; z) = \wt{\omega}(x; \chi; s) = \dint{K} \chi(f(k\cdot x)) \abs{f(k\cdot x)}^{s+\ve} dk,
$$
where $\chi$ is a character of $k^\times/k^{\times 2}$, $\ve = (-\frac12, \ldots, -\frac12)$, and $z$ is related to $s$ by
$$
\begin{array}{ll}
z_1 = s_1 + s_2 + s_3 + s_4, & z_2 = s_3 + s_4,\\
z_3 = s_1 + s_3, & z_4 = s_2 + s_3.
\end{array}
$$
$W = W_0 \cong (S_2 \ltimes (C_2)^2) \times (S_2)^2$, and (A5) is satisfied for each simple root.
Calculating functional equations by using zeta integrals of type $\wt{\Omega}(\wt{\phi}, \chi; s)$ (in Remark~3.15), we have explicit formulas of spherical functions on $X$. 
By spherical transform we see $\calS(K\backslash X)$ is a free $\hec$-module of rank $4$. 
Each spherical function is associated to some $z \in \C^4/W$ through $\wt{\lam}_z$, and the space of spherical functions associated to $\wt{\lam}_z$ has dimension $4$ and there is a basis explicitly given by terms of $\omega(x;\chi;z)$.
The small prehomogeneous vector spaces are isomorphic to $(GL_1 \times GL_1, V)$ over $k$ (cf. \cite{Hamb}-\S 4.1), and those functional equations are reduced to Tate's formula(\cite{Tate}-\S 2).

\vspace{2cm}
%\newpage
%\input{refer-hamb}

\bibliographystyle{amsalpha}

\vspace{15mm}
\begin{flushleft}
Yumiko HIRONAKA\\

\vspace{5mm}
Department of Mathematics\\ 
Faculty of Education and Integrated Sciences\\
Waseda University\\
Nishi-Waseda, Tokyo 169-8050, JAPAN

\makeatletter

\vspace{3mm}
e-mail : hironaka@waseda.jp
\end{flushleft}

\end{document}